\documentclass{article}
\usepackage[english]{babel}

\usepackage{amssymb}
\usepackage{amsmath}
\usepackage{amsthm}
\usepackage{mathcomp}
\usepackage{graphicx}
\usepackage{paralist} 
\usepackage{url}
\usepackage{caption}

\usepackage{enumitem}
\usepackage{color}

\usepackage{setspace}
\title{Universality of Composition Operators and Applications to Holomorphic Dynamics}
\AtBeginDocument{}

\theoremstyle{plain}
\newtheorem{theo}{Theorem}[section]
\newtheorem{lem}[theo]{Lemma}
\newtheorem{cor}[theo]{Corollary}

\theoremstyle{definition}
\newtheorem{defi}[theo]{Definition}

\newtheorem{rem-nota}[theo]{Remark and Notation}
\newtheorem{rem-defi}[theo]{Remark and Definition}

\newcommand{\dist}{\mbox{dist}}
\newcommand{\wt}{\widetilde}
\newcommand{\wh}{\widehat}
\newcommand{\ol}{\overline}
\newcommand{\bs}{\backslash}
\newcommand{\id}{\mbox{id}}

\newcommand{\bull}{\,\begin{picture}(-1,1)(-1,-3)\circle*{3}\end{picture}\,\,\,}
\newcommand{\Chi}{\raisebox{2pt}{$\chi$}}
\allowdisplaybreaks

\begin{document}



\begin{center}
\Large\textbf{Universality of Composition Operators and Applications to Holomorphic Dynamics}
\vspace{5pt}
\end{center}

\begin{center}
\large Andreas Jung\footnote{The author has been supported by the Stipendienstiftung Rheinland-Pfalz. Address:\\ University of Trier, FB IV, Mathematics, 54286 Trier, Germany, s4anjung@uni-trier.de}
\end{center}

\vspace{5pt}


\begingroup
\leftskip1.2cm
\rightskip\leftskip
\small
\textbf{Abstract.} 
By investigating which level of universality composition operators $C_f$ can have, where the symbol $f$ is given by the restriction of an entire function to suitable parts of the Fatou set of $f$, this work combines the theory of dynamics of continuous linear operators on spaces of holomorphic functions with the theory of non-linear complex dynamics on the complex plane.
\par
\endgroup

\vspace{10pt}



\section{Introduction}
\label{sec:Introduction}

\normalsize

For topological spaces $X$, $Y$ and a family $\{T_\iota:\iota\in I\}$ of continuous functions $T_\iota:X\rightarrow Y$, an element $x\in X$ is called \textit{universal} for $\{T_\iota:\iota\in I\}$ if the set
$\{T_\iota(x):\iota\in I\}$
is dense in $Y$. 
In case that $X$ is a topological vector space and $T:X\rightarrow X$ is a continuous linear operator, a vector $x\in X$ is called \textit{hypercyclic} for $T$ if it is universal for the family $\{T^n:n\in\mathbb{N}\}$, i.e.\,\,if its orbit 
$\{T^n(x):n\in\mathbb{N}\}$
is dense in $X$ (where for sets $X$ and functions $F:X\rightarrow X$, $F^n:=F\circ\ldots\circ F$ denotes the $n$-th iterate of $F$, that is $F^0:=\id_X$ and $F^n:=F\circ F^{n-1}$ for $n\in\mathbb{N}$). We investigate which level of universality sequences $(C_{f_n})$ can have, where 
$$C_f:H(D)\rightarrow H(D),\,\,\,C_f(g):=g\circ f,$$
denotes the \textit{composition operator} with \textit{symbol} $f$ for open sets $D\subset\mathbb{C}$ and holomorphic self-maps $f$ of $D$.
As usual, the space $H(D)$ of complex-valued holomorphic functions on $D$ is endowed with the topology of locally uniform convergence. 
The transform $C_f$ is a continuous linear operator on the Fr\'{e}chet space $H(D)$ and we have $(C_f)^n=C_{f^n}$ for all $n\in\mathbb{N}$.\medskip

The first work on universality of composition operators on general open subsets of the complex plane is due to W.\,\,Luh in 1993 (see \cite{luh}, Theorem on p.\,161). In case that all components of $D$ are simply connected, he proved the existence of a sequence $(f_n)$ of conformal self-maps of $D$ and the existence of a function $g\in H(D)$ having the property that for each compact subset $K$ of $D$ with connected complement the set $\{C_{f_n}(g)|_K:n\in\mathbb{N}\}$ is dense in the Banach space $\big(A(K),\|\bull\|_K\big)$, where $A(K)$ denotes the space of complex-valued continuous functions on $D$ which are holomorphic on $D^\circ$, and $\|\bull\|_K$ is the uniform norm on $K$, i.e.
$$\|h\|_K:=\max_{z\in K}|h(z)|,\,\,\,h\in A(K).$$
In 1995, L.\,\,Bernal-Gonz\'{a}lez and A.\,\,Montes-Rodr\'{i}guez were the first to characterize universality of a sequence of composition operators $(C_{f_n})$ for a given sequence $(f_n)$ of conformal self-maps of $D$ (cf.\,\cite{bernal-montes}). In order to state their main result, we introduce the following definition, which was coined by them:
\begin{defi}\label{run-away}
A sequence $(f_n)$ of holomorphic self-maps of a domain $\Omega\subset\mathbb{C}$ is called a \textit{run-away sequence} if for each compact subset $K$ of $\Omega$ there exists an $N\in\mathbb{N}$ such that $f_N(K)\cap K=\emptyset$.
\end{defi}
The main result of Bernal and Montes reads as follows (see \cite{bernal-montes}, Theorem 3.6): 
\begin{theo}\label{theorem_bernal_montes}\textnormal{(Bernal, Montes, 1995)}
Let $\Omega\subset\mathbb{C}$ be a domain and let $(f_n)$ be a sequence of conformal self-maps of $\Omega$.
\begin{itemize}
\itemsep0pt
\item[i)]If $\Omega$ is not conformally equivalent to $\mathbb{C}\bs\{0\}$, 
the sequence of composition operators $(C_{f_n})$ is universal if and only if $(f_n)$ is a run-away sequence. In case that one of these conditions holds, the set of universal functions for $(C_{f_n})$ is comeager in $H(\Omega)$.
\item[ii)]There exists a function $g\in H(\Omega)$ having the property that for each compact subset $K$ \hspace{-1pt}of \,$\Omega$ with connected complement the set $\left\{C_{f_n}(g)|_K:n\in\mathbb{N}\right\}$ is dense in $A(K)$ if and only if $(f_n)$ is a run-away sequence. In case of existence, the set of such functions is comeager in $H(\Omega)$.
\end{itemize}
\end{theo}
A subset of a Baire space is comeager if it contains a dense $G_\delta$-set. Bernal and Montes showed that the statement in i) does not hold for $\Omega=\mathbb{C}\bs\{0\}$ (see \cite{bernal-montes}, Remark on p.\,55). In 2009, K.-G.\,\,Grosse-Erdmann and R.\,\,Mortini characterized universality of composition operators without the assumption of conformality of the corresponding symbols -- they only required $(f_n)$ to be a sequence of holomorphic self-maps of $\Omega$ (cf.\,\cite{grosse-erdmann-mortini}). Their main result concerning the case of iteration of one holomorphic function is the following (cf.\,\cite{grosse-erdmann-mortini}, Theorem 3.21(a)):
\begin{theo}\label{theorem_grosse_erdmann_mortini}\textnormal{(Grosse-Erdmann, Mortini, 2009)}
Let $\Omega\subset\mathbb{C}$ be a simply connected domain and let $f$ be a holomorphic self-map of $\Omega$. Then $C_f$ is hypercyclic if and only if $(f^n)$ is a run-away sequence and $f$ is injective. 
\end{theo}

If one of the two conditions in Theorem \ref{theorem_grosse_erdmann_mortini} holds, the set of hypercyclic functions for $C_f$ is comeager in $H(\Omega)$. It is easy to see that injectivity of the symbol is necessary in order to obtain hypercyclicity of the corresponding composition operator. Indeed, if $f$ is not injective, 
we can choose points $z_1,z_2\in\Omega$ with $z_1\neq z_2$ and $f(z_1)=f(z_2)$. 
Then, for $g\in H(\Omega)$, only functions $h\in H(\Omega)$ with $h(z_1)=h(z_2)$ can be pointwise limit functions of a subsequence of $(g\circ f^n)$.\medskip 

Returning to the general situation of an open set $D\subset\mathbb{C}$ and a holomorphic self-map $f$ of $D$, we recall that a function $g\in H(D)$ is called hypercyclic for $C_f$ if the set $\left\{g\circ f^n:n\in\mathbb{N}\right\}$ is dense in $H(D)$. This is equivalent to the fact that for each $h\in H(D)$ there exists a strictly increasing sequence $(n_k)$ in $\mathbb{N}$ such that $(g\circ f^{n_k})$ converges to $h$ locally uniformly on $D$, i.e.\,\,each function in $H(D)$ can be approximated locally uniformly on $D$ via subsequences of $(g\circ f^n)$. For many situations which will occur in this work, it is useful to modify this definition of hypercyclicity with regard to the following two aspects:
\begin{itemize}
\itemsep0pt
\item[$\bullet$] We want to obtain hypercyclic functions for $C_f$ which are defined (and holomorphic) on larger open sets than $D$.
\item[$\bullet$] For hypercyclic functions $g$ for $C_f$ and suitable sets $M\subset D$, we want to approximate all elements of certain subclasses of the space $C(M$) of continuous functions on $M$ via subsequences of $\left(g\circ f^n|_M\right)$, i.e.\,\,we want to approximate classes of functions which are larger than $H(D)$ on sets which are smaller than $D$.
\end{itemize}

\vspace{2pt}
If $M\subset D$ is locally compact and $\sigma$-compact, there exists a compact exhaustion $(K_n)$ of $M$, and the family of seminorms $\{\|\bull\|_{K_n}:n\in\mathbb{N}\}$ turns $C(M)$ into a Fr\'{e}chet space. The induced topology does not depend on the choice of the sets $K_n$. In the following, the space $C(M)$ will always be endowed with this topology. Special cases are the spaces $A(K)$ if $M=K$ is compact and $H(U)$ if $M=U$ is open. In the latter case, a basis of the above topology is given by
$$\big\{V_{\varepsilon,L,U}(g):\varepsilon>0,L\subset U\mbox{ compact},g\in H(U)\big\},$$
where $V_{\varepsilon,L,U}(g):=\{h\in H(U):\|h-g\|_{L}<\varepsilon\}$. 
For $f_n,f\in C(M)$, $n\in\mathbb{N}$, the sequence $(f_n)$ converges to $f$ in the Fr\'{e}chet space $C(M)$ if and only if $(f_n)$ converges to $f$ locally uniformly on $M$. \medskip 

In order to comply with the two above-mentioned modifications, we now fix a domain $\Omega\subset\mathbb{C}$ with $\Omega\supset D$ and we introduce the following definition:
\begin{defi}\label{defi_F-universal}
Let $M\subset D$ be locally compact and $\sigma$-compact and let $\mathcal{F}\subset C(M)$ be a family of functions in the Fr\'{e}chet space $C(M)$. A function $g\in H(\Omega)$ is called \textit{$\Omega$-$\mathcal{F}$-universal} for $C_f$ if
$$\mathcal{F}\subset\ol{\left\{g\circ f^n|_M:n\in\mathbb{N}\right\}}^{C(M)}.$$
The composition operator $C_f$ is called \textit{$\Omega$-$\mathcal{F}$-universal} if there exists an $\Omega$-$\mathcal{F}$-universal function for $C_f$. If it is clear from the context which domain $\Omega$ is considered, $\Omega$-$\mathcal{F}$-universality will briefly be referred to as $\mathcal{F}$-universality.
\end{defi}

In case of $\left\{g\circ f^n|_M:n\in\mathbb{N}\right\}\subset\mathcal{F}$, the condition in Definition \ref{defi_F-universal} is equivalent to the fact that the set $\left\{g\circ f^n|_M:n\in\mathbb{N}\right\}$ is dense in $\mathcal{F}$.
For $\Omega=D$, a function $g\in H(D)$ is hypercyclic for $C_f$ if and only if it is $H(D)$-universal for $C_f$.\medskip

In general, we are mainly interested in the following question: For open subsets $U$ of $D$ (compact subsets $K$ of $D$), we want to find conditions on $f$ and $U$ (on $K$) such that the composition operator $C_f$ is $H(U)$-universal ($A(K)$-universal). In Section \ref{sec:General Universality Statements for Composition Operators with Holomorphic Symbol}, we will formulate and prove several statements of this kind. The main motivation of this work consists in applying such general universality statements to situations which naturally occur in the theory of complex dynamics, where the long-time behaviour of certain discrete topological dynamical systems is investigated. The systems which are studied are given by the iteration of a transcendental entire function or a rational function. To simplify matters, we will only consider entire functions in this work -- however, all following results also hold for rational functions (cf.\,\cite{jung}). The local behaviour of the considered function near fixed points plays an important role in the study of the long-time behaviour of the corresponding dynamical system. For a complex-valued function $g$ which is holomorphic near a point $z_0\in\mathbb{C}$ with $g(z_0)=z_0$, the \textit{multiplier} of $g$ at $z_0$ is defined as $\lambda:=g'(z_0)$, and the fixed point $z_0$ is called \textit{superattracting}, \textit{attracting}, \textit{indifferent} or \textit{repelling} if $\lambda=0$, $0<|\lambda|<1$, $|\lambda|=1$ or $|\lambda|>1$, respectively. Indifferent fixed points with $\lambda=1$ are also called \textit{neutral}.\medskip


We denote by $\mathbb{C}_\infty:=\mathbb{C}\cup\{\infty\}$ the extended complex plane and we endow this set with the chordal metric $\Chi$. For a topological space $X$ and a metric space $Y$, a family $\mathcal{F}\subset C(X,Y)$ is called \textit{normal} if each sequence in $\mathcal{F}$ has a subsequence which converges uniformly on all compact subsets of X to a function in the space $C(X,Y)$ of continuous functions from $X$ to $Y$. If $X$ is locally compact and has a compact exhaustion, normality is a local property, i.e.\,\,a family $\mathcal{F}\subset C(X,Y)$ is normal if and only if for each point $x\in X$ there exists a neighbourhood $W$ of $x$ such that $\big\{f|_{W}:f\in\mathcal{F}\big\}$ is a normal family in $C(W,Y)$ (cf.\,\cite{schiff}, Theorem 2.1.2, in case that $X$ is a domain in $\mathbb{C}$).\medskip

For the rest of this section, we now fix an entire function $f$ which is not a polynomial of degree $0$ or $1$. The following definition is crucial in the theory of complex dynamics:
\begin{defi}
The \textit{Fatou set} $F_f$ of $f$ is defined as the set of all points $z\in\mathbb{C}$ for which there exists a neighbourhood $W$ of $z$ such that $\left\{f^n|_W:n\in\mathbb{N}\right\}$ is a normal family in $C(W,\mathbb{C}_\infty)$. The complement of $F_f$ is called the \textit{Julia set} of $f$ and it is denoted by $J_f$.
\end{defi}

By definition, $F_f$ is an open set. Thus, $J_f$ is closed. Due to the above considerations, $\{f^n|_{F_f}:n\in\mathbb{N}\}$ is a normal family in $C(F_f,\mathbb{C}_\infty)$. It is well-known that the Julia set of $f$ is completely invariant under $f$ (see \cite{schleicher} p.\,299 and Theorem 1.7). Hence, the same is also true for the Fatou set of $f$. Moreover, $J_f$ is a perfect set, i.e.\,\,$J_f$ has no isolated points (see \cite{schleicher}, Theorem 1.4), and we have $F_{f^n}=F_f$ as well as $J_{f^n}=J_f$ for all $n\in\mathbb{N}$ (see \cite{schleicher}, p.\,299). (Super-)attracting fixed points of $f$ are contained in the Fatou set of $f$, whereas repelling and 
neutral fixed points of $f$ belong to the Julia set of $f$.\medskip


A reason for splitting the complex plane into the two disjoint subsets $F_f$ and $J_f$ is the following: It is well-known that the sequence of iterates $(f^n)$ behaves quite ``chaotically'' on $J_f$ whereas, on the other hand, the behaviour of $(f^n)$ on $F_f$ is quite ``nice'' and well-understood. The ``chaotic'' behaviour of $(f^n)$ on $J_f$ can be described as follows: For each point $z\in J_f$ and each arbitrarily small open neighbourhood $V$ of $z$, we have
$$\bigg|\mathbb{C}_\infty\bs\bigcup_{n\in\mathbb{N}}f^n(V)\bigg|\leq2.$$
Thus, every value in the extended complex plane, with at most two exceptions, is assumed on $V$ under the iteration of $f$. This result is a direct application of Montel's theorem 
which states that a family of meromorphic functions on a domain $G\subset\mathbb{C}$ omitting three fixed values in $\mathbb{C}_\infty$ is a normal family in $C(G,\mathbb{C}_\infty)$ (cf.\,\cite{carleson-gamelin}, Theorem I.3.2).\medskip

In order to point out that the sequence of iterates $(f^n)$ behaves ``nicely'' on $F_f$, we now consider a component $G$ of $F_f$. For each $n\in\mathbb{N}$, the continuity of $f^n$ and the invariance of $F_f$ under $f$ imply that $f^n(G)$ is a connected subset of $F_f$. Hence, there exists a component $G_n$ of $F_f$ with $f^n(G)\subset G_n$. Using the notation $G_0:=G$, we call $G$ \textit{preperiodic} if there exist integers $p>m\geq0$ with $G_p=G_m$. If $G$ is preperiodic and if in particular $m=0$, we have $f^p(G)\subset G$ and $G$ is called \textit{periodic} with \textit{period} $p$. Finally, if $G$ is periodic with period $1$, we have $f(G)\subset G$ and we call $G$ \textit{invariant}. 
The following theorem gives a complete description of the behaviour of the sequence of iterates $(f^n)$ on components of the Fatou set of $f$. It is mainly due to P.\,\,Fatou in 1919 and H.\,\,Cremer in 1932 (cf.\,\cite{bergweiler_uebersicht}, p.\,163).
\begin{theo}\label{classifcation_theorem}\textnormal{(Classification Theorem of Fatou Components)}
Let $G$ be an invariant component of $F_f$. Then $G$ is of exactly one of the following types:
\begin{itemize}
\itemsep0pt
\item[$\bullet$] $G$ contains a superattracting fixed point $z_0$ of $f$ and $f^n|_G\rightarrow z_0$ (\textit{B\"ottcher domain}),
\item[$\bullet$] $G$ contains an attracting fixed point $z_0$ of $f$ and $f^n|_G\rightarrow z_0$ (\textit{Schr\"oder domain}),
\item[$\bullet$] $\partial G$ contains a neutral fixed point $z_0$ of $f$ and $f^n|_G\rightarrow z_0$ (\textit{Leau domain}),
\item[$\bullet$] $f$ is conjugated on $G$ to an irrational rotation on $\mathbb{D}:=\{z\in\mathbb{C}:|z|<1\}$ (\textit{Siegel disc}), 
\item[$\bullet$] $f$ is transcendental entire and $f^n|_G\rightarrow\infty$ (\textit{Baker domain}).
\end{itemize}
\end{theo}
(Pre-)periodic components of $F_f$ are classified in exactly the same way. A component of $F_f$ which is not preperiodic is called a \textit{wandering domain}. It is well-known that polynomials do not have any wandering domains (see \cite{sullivan}, Theorem 1). In view of the classification theorem, there exist only a few possibilities of how the iterates of $f$ can behave on the Fatou set of $f$. For this reason, we stated above that the behaviour of $(f^n)$ on $F_f$ is ``nice'' and well-understood.\medskip

The main motivation of this work is the following: Composing the iterates from the left with one other suitable holomorphic function, we will see that this nice behaviour of the sequence of iterates on the Fatou set can ``change dramatically'' such that it becomes quite similarly chaotic as on the Julia set! The general aim of this work is to study the long-time behaviour of such modified sequences on the Fatou set. More precisely, this means that we will prove $H(U)$- and $A(K)$-universality of the composition operator $C_f$ for suitable open subsets $U$ of $F_f$ and suitable compact subsets $K$ of $F_f$. This will be the ``dramatic change of behaviour'' of the sequence of iterates $(f^n)$ composed with some universal function $g$ for $C_f$ -- because, in this situation, the family of compositions $\{g\circ f^n:n\in\mathbb{N}\}$ will be dense in some space of holomorphic functions, which means that the sequence of compositions $(g\circ f^n)$ will be kind of ``maximally divergent''. 
But this stands in sharp contrast to the fact that the sequence of iterates $(f^n)$ itself converges, behaves like a rotation or is wandering on all components of $F_f$, as stated in the classification theorem.\medskip

Moreover, if $U\subset F_f$ is a domain such that there exist a compact set $K\subset U$ which has an accumulation point in $U$ and a function $g$ which is $A(K)$-universal for $C_f$, the family of compositions $\{g\circ f^n|_U:n\in\mathbb{N}\}$ is not normal in $C(U,\mathbb{C}_\infty)$ any more (by Mergelyan's theorem, this is in particular fulfilled if $C_f$ is $H(U)$-universal). Indeed, choosing some $h\in A(K)$ which can not be extended holomorphically to $U$, the denseness of $\{g\circ f^n|_K:n\in\mathbb{N}\}$ in $A(K)$ yields a sequence $(n_j)$ in $\mathbb{N}$ with uniform convergence $g\circ f^{n_j}\rightarrow h$ on $K$. Assuming that the family $\{g\circ f^n|_U:n\in\mathbb{N}\}$ is normal in $C(U,\mathbb{C}_\infty)$, the Vitali-Porter theorem would imply the existence of a function $\psi\in H(U)\cup\{\infty\}$ with locally uniform convergence $g\circ f^{n_j}\rightarrow\psi$ on $U$, and we would obtain $\psi|_K=h$, a contradiction. 


\medskip


\section{General Universality Statements for Composition Operators with Holomorphic Symbol}
\label{sec:General Universality Statements for Composition Operators with Holomorphic Symbol}

Throughout this section, we fix an open set $D\subset\mathbb{C}$, a domain $\Omega\subset\mathbb{C}$ with $\Omega\supset D$ and a holomorphic function $f:D\rightarrow D$. As already mentioned in the introduction, our aim will be the following: For open subsets $U$ of $D$ and compact subsets $K$ of $D$, we want to find conditions on $f$ and $U$ (and on $K$, respectively) such that the composition operator $C_f$ is $H(U)$-universal ($A(K)$-universal).\medskip

As seen above, topological properties of sets play an important role for obtaining universality of a given sequence of composition operators, i.e.\,\,the above-mentioned theorems of Luh, Bernal and Montes as well as Grosse-Erdmann and Mortini hold for sets which have connected complement.
For a set $M\subset\mathbb{C}$, we call a component $H$ of $\mathbb{C}_\infty\bs M$ with $\infty\notin H$ a \textit{hole} of $M$. Moreover, we denote
\begin{eqnarray*}
\mathcal{K}(M)&:=&\left\{K\subset M:K\mbox{ compact, }K\neq\emptyset\right\},\\
\mathcal{K}_0(M)&:=&\left\{K\in\mathcal{K}(M):K\mbox{ has no holes}\right\},\\
\mathcal{U}(M)&:=&\left\{U\subset M:U\mbox{ open in }\mathbb{C},\,U\neq\emptyset\right\},\\
\mathcal{U}_0(M)&:=&\left\{U\in\mathcal{U}(M):U\mbox{ has no holes}\right\}.
\end{eqnarray*}

We remark that a subset $M$ of the complex plane has no holes if and only if $\mathbb{C}_\infty\bs M$ is connected. In this work, we will often use the fact that for an open set $U\subset\mathbb{C}$ and an injective holomorphic function $\varphi:U\rightarrow\mathbb{C}$ the image $\varphi(U)$ has the same number of holes as $U$. An analogous statement holds if $K\subset\mathbb{C}$ is compact and $\varphi$ is injective and holomorphic on an open neighbourhood of $K$, i.e.\,\,then $K$ and $\varphi(K)$ have the same number of holes (see \cite{remmert2}, p.\ 276 and cf.\,\cite{grosse-erdmann-mortini}, p.\ 360). Whenever we use one of these statements, we call this the \textit{invariance of the number of holes}. Sometimes, we want to ``fill up'' holes of compact sets relative to a given open superset. In order to specify what this means, we denote for $U\in\mathcal{U}(\mathbb{C})$ and $K\in\mathcal{K}(U)$ by
$\wh{K}_U:=\{z\in U:|f(z)|\leq||f||_K\mbox{ for all }f\in H(U)\}$
the \textit{holomorphically convex hull} of $K$ with respect to $U$, and we say that $K$ is \textit{$U$-convex} if $K=\wh{K}_U$. By definition, 
$\wh{K}_U$ is $U$-convex. It is well-known that the compact set $\wh{K}_U$ is the union of $K$ with all holes of $K$ which are contained in $U$ (cf.\,\cite{remmert2}, p.\,263f.). Thus, $K$ is $U$-convex if and only if each hole of $K$ contains a point in $\mathbb{C}\bs U$. If, in particular, $U=\mathbb{C}$, we write $\wh{K}:=\wh{K}_\mathbb{C}$ for the \textit{polynomially convex hull} of $K$ (cf.\,\cite{remmert2}, p.\,266).\medskip

In view of Theorem \ref{theorem_grosse_erdmann_mortini}, the following definition is now reasonable:
\begin{defi}\label{U_1}
We denote by $\mathcal{U}_0(D,\Omega,f)$ the set of all sets $U\in\mathcal{U}_0(D)$ for which we have 
locally uniform convergence $f^n|_U\rightarrow\partial_\infty\Omega$ (i.e.\,\,the sequence $(\dist(f^n(\bull),\partial_\infty\Omega))_{n\in\mathbb{N}}$ converges to $0$ locally uniformly on $U$) and injectivity of each iterate $f^n$ on $U$.
\end{defi}

Here, $\partial_\infty$ denotes the boundary of subsets of $\mathbb{C}$ with respect to $(\mathbb{C}_\infty,\Chi)$. 
The first condition in Definition \ref{U_1} corresponds to the run-away property of the sequence of iterates $(f^n)$ as defined by Bernal and Montes (see Definition \ref{run-away}) and as stated in Theorem \ref{theorem_grosse_erdmann_mortini}. At first view, the second condition seems to be stronger than the condition of injectivity of the symbol of the composition operator in Grosse-Erdmann's and Mortini's theorem. However, in the situation of Theorem \ref{theorem_grosse_erdmann_mortini}, the symbol is a self-map of the considered domain so that its injectivity implies the injectivity of all of its iterates, which is crucial for the proof of Theorem \ref{theorem_grosse_erdmann_mortini}.\medskip

Now, we can state and prove the following universality result, which will be the basis for several further universality statements. The proof will run along the lines of the proof of Theorem 3.2, (c)$\Longrightarrow$(a), in \cite{grosse-erdmann-mortini}.
\begin{theo}\label{central_universality_theorem}
Let $U\in\mathcal{U}_0(D,\Omega,f)$. Then the set of all functions in $H(\Omega)$ which are $H(U)$-universal for $C_f$ is a comeager set in $H(\Omega)$.
\end{theo}
\textbf{Proof:} 
For $n\in\mathbb{N}$, we consider the composition operator
$$C_{f^n,U}:H(\Omega)\rightarrow H(U),\,C_{f^n,U}(g):=g\circ f^n|_U.$$
According to the universality criterion, it suffices to show that the sequence $(C_{f^n,U})$ is topologically transitive (cf.\,\cite{grosse-erdmann-survey}, Theorem 1). In order to do so, let $\emptyset\neq\mathcal{V}\subset H(\Omega)$ and $\emptyset\neq\mathcal{W}\subset H(U)$ be open. 
For $g\in\mathcal{V}$, there exist $\varepsilon_1>0$ and $K_1\in\mathcal{K}(\Omega)$ such that $V_{\varepsilon_1,K_1,\Omega}(g)\subset\mathcal{V}$. Defining $$K:=\wh{(K_1)}_\Omega,$$ we obtain that $K$ is $\Omega$-convex with $V_{\varepsilon_1,K,\Omega}(g)=V_{\varepsilon_1,K_1,\Omega,}(g)\subset\mathcal{V}.$ 
For $h\in\mathcal{W}$, there exist $\varepsilon_2>0$ and $L_1\in\mathcal{K}(U)$ with $V_{\varepsilon_2,L_1,U}(h)\subset\mathcal{W}$. Considering the polynomially convex hull $$L:=\widehat{L_1},$$ we have
$V_{\varepsilon_2,L,U}(h)=V_{\varepsilon_2,L_1,U}(h)\subset\mathcal{W}.$
As $U$ has no holes, it follows that $L\in\mathcal{K}_0(U)$. Putting \mbox{$\varepsilon:=\min(\varepsilon_1,\varepsilon_2)>0$}, we obtain
$$V_{\varepsilon,K,\Omega}(g)\subset V_{\varepsilon_1,K,\Omega}(g)\subset\mathcal{V}\quad\mbox{ and }\quad V_{\varepsilon,L,U}(h)\subset V_{\varepsilon_2,L,U}(h)\subset\mathcal{W}.$$
Now, let $\delta:=\dist(K,\partial_\infty\Omega)>0$. Because of $U\in\mathcal{U}_0(D,\Omega,f)$, we obtain the uniform convergence $f^n|_L\rightarrow\partial_\infty\Omega$. Thus, there exists an $N\in\mathbb{N}$ with $\dist(f^N(z),\partial_\infty\Omega)<\delta$ for all $z\in L$, which implies $K\cap f^N(L)=\emptyset.$
According to $U\in\mathcal{U}_0(D,\Omega,f)$, the restriction $f^N|_U$ is injective. Hence, as $L$ has no holes, the invariance of the number of holes implies that $f^N(L)$ has no holes as well. Therefore, the disjoint union of the compact set $K$ with the compact set $f^N(L)$ does not produce a new hole so that the $\Omega$-convexity of $K$ implies that $K\cup f^N(L)$ is $\Omega$-convex again. Thus, we can choose from each hole of $K\cup f^N(L)$ a point which lies in $\mathbb{C}\bs\Omega$. Let $A$ be the union of these points and let $B:=A\cup\{\infty\}$. Then we have $B\subset\mathbb{C}_\infty\bs\Omega\subset\mathbb{C}_\infty\bs(K\cup f^N(L))$ and $B\cap C\neq\emptyset$ for all components $C$ of $\mathbb{C}_\infty\bs(K\cup f^N(L))$. We now consider the function
$$\varphi:K\cup f^N(L)\rightarrow\mathbb{C},\,\,\,\varphi(z):=\left\{\begin{array}{ll} g(z), & \mbox{if }z\in K\\ h\big((f^N|_U)^{-1}(z)\big), & \mbox{if }z\in f^N(L)\end{array}\right..$$
As the disjoint sets $K$ and $f^N(L)$ are compact with $K\subset\Omega$ and $f^N(L)\subset f^N(U)$ and since we have $g\in H(\Omega)$ as well as $h\in H(U)$, we see that $\varphi$ can be extended holomorphically to an open neighbourhood of $K\cup f^N(L)$. Hence, Runge's theorem yields a rational function $R$ having poles only in $B\subset\mathbb{C}_\infty\bs\Omega$ such that
$\|\varphi-R\|_{K\cup f^N(L)}<\varepsilon.$
The restriction $R|_\Omega$ is holomorphic on $\Omega$, and due to $\varphi=g$ on $K$ we obtain
$$\|g-R\|_K=\|\varphi-R\|_K\leq\|\varphi-R\|_{K\cup f^N(L)}<\varepsilon$$
and thus $R|_\Omega\in V_{\varepsilon,K,\Omega}(g)\subset\mathcal{V}$. Because of $\varphi=h\circ(f^N|_U)^{-1}$ on $f^N(L)$, we further obtain
\begin{eqnarray*}
\left\|h-C_{f^N,U}\left(R|_\Omega\right)\right\|_L&=&\left\|h-R\circ f^N\right\|_L=\left\|h\circ(f^N|_U)^{-1}-R\right\|_{f^N(L)}\\
&=&\|\varphi-R\|_{f^N(L)}\leq\|\varphi-R\|_{K\cup f^N(L)}<\varepsilon
\end{eqnarray*}
and thus $C_{f^N,U}\left(R|_\Omega\right)\in V_{\varepsilon,L,U}(h)\subset\mathcal{W}$. Altogether, it follows that $C_{f^N,U}\left(R|_\Omega\right)$ is contained in the intersection $C_{f^N,U}(\mathcal{V})\cap\mathcal{W}$
so that the topological transitivity of the sequence $(C_{f^n,U})$ is shown.\hfill$\square$\medskip

The difference between the above theorem and the result of Grosse-Erdmann and Mortini (Theorem \ref{theorem_grosse_erdmann_mortini}) is the following: Applying Theorem \ref{central_universality_theorem}, we obtain $H(U)$-universal functions for $C_f$ which are actually defined and holomorphic on the ``large'' domain $\Omega$ -- and not only on the set $U$ itself which is contained in the open set $D$ on which the symbol $f$ is defined. Now, it is our aim to extend the statement of Theorem \ref{central_universality_theorem} in the following way: We want to find conditions such that the composition operator $C_f$ is $H(U)$-universal for all sets $U\in\mathcal{U}_0(D)$. This intention originates from the above-mentioned second result of Bernal and Montes (Theorem \ref{theorem_bernal_montes}\,ii)). In view of Theorem \ref{central_universality_theorem}, it is natural to require the existence of a sequence in $\,\mathcal{U}_0(D,\Omega,f)$, which fulfils the following property in $D$:
\begin{defi}
Let $V\in\mathcal{U}(\mathbb{C})$. A sequence $(M_n)$ of subsets of $V$ is called \textit{$\mathcal{K}_0(V)$-exhausting} if for each set $L\in\mathcal{K}_0(V)$ there exists an $N\in\mathbb{N}$ with $L\subset M_N^{\circ}$.
\end{defi}

Each open set $V\subset\mathbb{C}$ has a $\mathcal{K}_0(V)$-exhausting sequence in $\,\mathcal{U}_0(V)$. Indeed, Bernal and Montes have shown that each domain $G\subset\mathbb{C}$ has a $\mathcal{K}_0(G)$-exhausting sequence in $\mathcal{K}_0(G)$ (see \cite{bernal-montes}, Lemma 2.9). Considering their proof of this statement, one can verify that it also holds for arbitrary open subsets of the complex plane. Hence, there exists a $\mathcal{K}_0(V)$-exhausting sequence $(K_n)$ in $\mathcal{K}_0(V)$. As we can find for each $n\in\mathbb{N}$ a set $U_n\in\mathcal{U}_0(V)$ with $U_n\supset K_n$ (see e.g.\,\cite{grosse-erdmann-diss}, Satz 2.2.4 and the first part of its proof), we obtain that $(U_n)$ is a $\mathcal{K}_0(V)$-exhausting sequence in $\,\mathcal{U}_0(V)$. 
\begin{cor}\label{central_universality_theorem_all_open}
Let there exist a $\mathcal{K}_0(D)$-exhausting sequence in $\,\mathcal{U}_0(D,\Omega,f)$. Then the set of all functions in $H(\Omega)$ which are $H(U)$-universal for $C_f$ for all $U\in\mathcal{U}_0(D)$ is a comeager set in $H(\Omega)$.
\end{cor}
\textbf{Proof:} 
Let $(U_m)$ be a $\mathcal{K}_0(D)$-exhausting sequence in $\,\mathcal{U}_0(D,\Omega,f)$. Moreover, for $m\in\mathbb{N}$, let $\mathcal{G}_m$ be the set of all functions in $H(\Omega)$ which are $H(U_m)$-universal for $C_f$. Due to Theorem \ref{central_universality_theorem}, each set $\mathcal{G}_m$ is comeager in $H(\Omega)$ so that the same is also true for the countable intersection
$\mathcal{G}:=\bigcap_{m\in\mathbb{N}}\mathcal{G}_m.$
Let $\mathcal{H}$ be the set of all functions in $H(\Omega)$ which are $H(U)$-universal for $C_f$ for all $U\in\mathcal{U}_0(D)$. We now show the inclusion $\mathcal{G}\subset\mathcal{H}$, which will complete the proof. For this purpose, let $g\in\mathcal{G}$ and $U\in\mathcal{U}_0(D)$. As we have to show the denseness of the set $\left\{g\circ f^n|_U:n\in\mathbb{N}\right\}$ in $H(U)$, let moreover $h\in H(U), K\in\mathcal{K}(U)$ and $\varepsilon>0$. Because of $U\in\mathcal{U}_0(D)\subset\mathcal{U}_0(\mathbb{C})$, we have $\widehat{K}\in\mathcal{K}_0(U)\subset\mathcal{K}_0(D)$. Runge's theorem implies the existence of a polynomial $P$ with $\|P-h\|_{\widehat{K}}<\varepsilon\slash2$. As $(U_m)$ is $\mathcal{K}_0(D)$-exhausting, there exists an $M\in\mathbb{N}$ with $\widehat{K}\subset U_{M}$. Due to $g\in\mathcal{G}\subset\mathcal{G}_{M}$ and $P|_{U_{M}}\in H(U_{M})$, we obtain an $N\in\mathbb{N}$ with
$\big\|g\circ f^N|_{U_{M}}-P\big\|_{\widehat{K}}<\varepsilon\slash2.$
Altogether, it follows
$$\left\|g\circ f^N|_U-h\right\|_K\leq\big\|g\circ f^N|_{U_{M}}-P\big\|_{\widehat{K}}+\|P-h\|_{\widehat{K}}<\frac{\varepsilon}{2}+\frac{\varepsilon}{2}=\varepsilon.$$
Thus, we have shown that $\left\{g\circ f^n|_U:n\in\mathbb{N}\right\}$ is dense in $H(U)$, which means that $g$ is $H(U)$-universal for $C_f$. Therefore, we obtain $g\in\mathcal{H}$ so that $\mathcal{G}\subset\mathcal{H}$ is shown.\hfill$\square$\medskip

In particular, the assumptions of Corollary \ref{central_universality_theorem_all_open} are fulfilled if $f$ is injective and if $f^n|_D\rightarrow\partial_\infty\Omega$ locally uniformly. Thus, we obtain:
\begin{cor}\label{central_universality_theorem_all_open_special}
Let $f$ be injective and let $f^n\rightarrow\partial_\infty\Omega$ locally uniformly on $D$. Then the set of all functions in $H(\Omega)$ which are $H(U)$-universal for $C_f$ for all $U\in\mathcal{U}_0(D)$ is a comeager set in $H(\Omega)$.
\end{cor}

In the special case that $D$ has no holes, we obtain the following result:
\begin{cor}\label{central_universality_theorem_all_open_special_simply_connected}
Let $D\in\mathcal{U}_0(\mathbb{C})$, let $f$ be injective and let $f^n\rightarrow\partial_\infty\Omega$ locally uniformly on $D$. Then the set of all functions in $H(\Omega)$ which are $H(D)$-universal for $C_f$ is a comeager set in $H(\Omega)$.
\end{cor}

Similarly to Definition \ref{U_1}, the following definition is reasonable in case of finite subsets of D:
\begin{defi}
We denote by $\mathcal{E}(D,\Omega,f)$ the set of all finite subsets $E$ of $D$ for which we have convergence $f^n|_E\rightarrow\partial_\infty\Omega$ and injectivity of each iterate $f^n$ on $E$.
\end{defi}

The same approach as in the proof of Theorem \ref{central_universality_theorem} now yields the following result:
\begin{cor}\label{central_universality_theorem_finite}
Let $E\in\mathcal{E}(D,\Omega,f)$. Then the set of all functions in $H(\Omega)$ which are $C(E)$-universal for $C_f$ is a comeager set in $H(\Omega)$.
\end{cor}

\section{Applications: Local Theory}    
\label{sec:Applications: Local Theory} 

In this section, we want to state some first applications of the general universality results of the previous section to the theory of complex dynamics. These applications will all be of local nature, meaning that we will prove universality of composition operators with locally defined symbols which are holomorphic near a fixed point. In view of the classification theorem of Fatou components, we will consider attracting, neutral and superattracting fixed points of the symbols -- nevertheless, the results in this section do not depend on results of the Fatou-Julia theory, i.e.\,\,here we consider symbols $f$ which only need to be defined locally near a point $z_0\in\mathbb{C}$ with $f(z_0)=z_0$.

\subsection{Attracting Fixed Points}
\label{subsec:Attracting Fixed Points}

Let $z_0$ be an attracting fixed point of $f$, i.e.\,\,for $\lambda:=f'(z_0)$ we have $0<|\lambda|<1$. The main consideration we need now is the concept of conformal conjugation. According to G.\,\,K\oe{}nigs' linearization theorem, there exist open neighbourhoods $U$ of $z_0$ and $V$ of $0$ as well as a conformal map $\varphi:U\rightarrow V$ which conjugates the map $f|_U:U\rightarrow U$ to the linear function $F:V\rightarrow V,\,F(w):=\lambda w$, i.e.
$$\varphi\circ f^n=F^n\circ\varphi=\lambda^n\cdot\varphi$$
holds on $U$ for all $n\in\mathbb{N}$ (see e.g.\,\cite{carleson-gamelin}, Theorem II.2.1). In particular, we obtain $\varphi(z_0)=0$ and hence locally uniform convergence $f^n\rightarrow z_0$ on $U$.
\begin{theo}\label{attracting_case}
In the above situation, the set of all functions in $H(\mathbb{C}\bs\{z_0\})$ which are $H(W)$-universal for $C_f$ for all $W\in\mathcal{U}_0(U\bs\{z_0\})$ is a comeager set in $H(\mathbb{C}\bs\{z_0\})$.
\end{theo}
\textbf{Proof:} 
We consider the domain $\Omega:=\mathbb{C}\bs\{z_0\}$ and the open set $D:=U\bs\{z_0\}\subset\Omega$. As $F,\varphi$ and $\varphi^{-1}$ are injective, the same is also true for \mbox{$f|_U=\varphi^{-1}\circ F\circ\varphi$.} Hence, the restriction $f|_D:D\rightarrow D$ is injective. Due to the locally uniform convergence $f^n|_D\rightarrow z_0\in\partial_\infty\Omega$, the assertion follows from Corollary \ref{central_universality_theorem_all_open_special}.\hfill$\square$\medskip

As the open set $U\bs\{z_0\}$ has a hole (namely $\{z_0\}$), one can show that there does not exist a function in $H(U\bs\{z_0\})$ which is $H(U\bs\{z_0\})$-universal for $C_f$ (cf.\,\,the Remark in \cite{bernal-montes} on p.\,\,55 and observe that $f|_{U\bs\{z_0\}}$ is conjugated to the linear function $w\mapsto\lambda w$ on $V\bs\{0\}$). 


\subsection{Neutral Fixed Points}
\label{subsec:Neutral Fixed Points}

Let $z_0$ be a neutral fixed point of $f$, i.e.\,\,we have $f'(z_0)=1$, and let $f$ not be the identity map. Putting $$m:=\min\{n\in\mathbb{N}:f^{(n+1)}(z_0)\neq0\},$$ 
an application of the Leau-Fatou flower theorem implies the existence of $m$ \textit{attracting petals} $P_1,\ldots,P_m$ for $f$, which are pairwise disjoint simply connected domains in $\mathbb{C}$ with $f(P_k)\subset P_k$, $z_0\in\partial P_k$ and $f^n|_{P_k}\rightarrow z_0$ uniformly, such that for each $k\in\{1,\ldots,m\}$, there exists a simply connected domain $V_k\subset\mathbb{C}$ containing a right half-plane 
as well as a conformal map $\varphi_k:P_k\rightarrow V_k$ which conjugates $f|_{P_k}$ to the translation $w\mapsto w+1$ on $V_k$, i.e.\,\,the equation
$$\varphi_k\circ f^n=\varphi_k+n$$
holds on $P_k$ for all integers $n\in\mathbb{N}$ (see e.g.\cite{milnor}, Theorems 10.5 and 10.7). Defining $P:=\bigcup_{k=1}^mP_k$ to be the union of all attracting petals, the following universality result holds:
\begin{theo}\label{neutral_case}
In the above situation, the set of all functions in $H(\mathbb{C}\bs\{z_0\})$ which are $H(P)$-universal for $C_f$ is a comeager set in $H(\mathbb{C}\bs\{z_0\})$.
\end{theo}
\textbf{Proof:} 
We consider the domain $\Omega:=\mathbb{C}\bs\{z_0\}$ and the open set $D:=P\in\mathcal{U}_0(\Omega)$. First, we show that $f|_D:D\rightarrow D$ is injective. To this end, let $z,w\in D$ with $f(z)=f(w)$. Then there exist $k,l\in\{1,\ldots,m\}$ with $z\in P_k$ and $w\in P_l$, and the invariance of $P_k$ and $P_l$ under $f$ implies $f(z)\in P_k\cap P_l$. As $P_k$ and $P_l$ are disjoint, we obtain $k=l$ so that the injectivity of $f|_{P_k}=\varphi_k^{-1}\circ(\varphi_k+1)$ yields $z=w$. Due to the uniform convergence $f^n|_D\rightarrow z_0\in\partial_ \infty\Omega$, the assertion follows from Corollary \ref{central_universality_theorem_all_open_special_simply_connected}.\hfill$\square$


\subsection{Superattracting Fixed Points}
\label{subsec:Superattracting Fixed Points}

Let $z_0$ be a superattracting fixed point of $f$, i.e.\,\,we have $f'(z_0)=0$, and let $f$ not be constant. Then, for $p:=\min\{n\in\mathbb{N}:f^{(n)}(z_0)\neq0\}$, we have $p\geq2$. 
Again, the concept of conformal conjugation will be the main consideration now. Due to B\"ottcher's theorem, there exist open neighbourhoods $U$ of $z_0$ and $V$ of $0$ as well as a conformal map $\varphi:U\rightarrow V$ which conjugates $f|_U:U\rightarrow U$ to the $p$-th monomial $Q$ on $V$, i.e.
$$\varphi(f^n(z))=Q^n(\varphi(z))=(\varphi(z))^{p^n}$$  
holds for all $z\in U$ and for all $n\in\mathbb{N}$ (see e.g.\,\cite{carleson-gamelin}, Theorem II.4.1). Without loss of generality, we may assume $V\subset\mathbb{D}$ so that the equality $\varphi(z_0)=\varphi(z_0)^p$ implies $\varphi(z_0)=0$. Thus, analogously to the considerations before Theorem \ref{attracting_case}, we obtain locally uniform convergence $f^n\rightarrow z_0$ on $U$.\medskip

Similarly to the previous two subsections, it is now our aim to formulate and prove a universality result for the composition operator $C_f$ which holds locally near $z_0$. The problem compared to the situations of attracting or neutral fixed points is the following: As we have $f(z)=\varphi^{-1}\big((\varphi(z))^p\big)$ for all $z\in U$, we see that the symbol $f$ now is not injective on $U$. Hence, we cannot apply Corollaries \ref{central_universality_theorem_all_open_special} or \ref{central_universality_theorem_all_open_special_simply_connected} for proving $H(W)$-universality of $C_f$ for suitable sets $W\in\mathcal{U}_0(U\bs\{z_0\})$. However, for $n\in\mathbb{N}$, the $p^n$-th monomial is injective on each angular sector of the complex plane with arc length smaller than $2\pi\slash p^n$. In view of Theorem \ref{central_universality_theorem} and the considerations after Definition \ref{U_1}, we need injectivity of all iterates of $f$ in order to obtain universality of $C_f$. 
Because of $2\pi\slash p^n\rightarrow0$ and due to the above conjugation, this means that for sets $M\subset U\bs\{z_0\}$ and families $\mathcal{F}\subset C(M)$, we only have a chance of obtaining $\mathcal{F}$-universality of $C_f$ if $M$ has empty interior (cf.\,\,the considerations after Theorem \ref{superattracting_1}). It will be our aim to prove that $C_f$ is $A(K)$-universal for suitable compact subsets $K$ of $U\bs\{z_0\}$ with $K^\circ=\emptyset$. As compensation for only considering such ``small'' subsets, we want to obtain ``many'' of them. In order to specify what this means, we consider for $\delta>0$ with $U_\delta[0]:=\{z\in\mathbb{C}:|z|\leq\delta\}\subset V$ the set
$B_\delta:=\varphi^{-1}(U_\delta[0]),$  
which is a compact neighbourhood of $z_0$ that is contained in $U$. Then we have $\varphi(B_\delta^\circ)=\varphi(B_\delta)^\circ=U_\delta(0):=\{z\in\mathbb{C}:|z|<\delta\}$ so that
$$\varphi|_{B_\delta^{\circ}\bs\{z_0\}}:B_\delta^{\circ}\bs\{z_0\}\rightarrow U_\delta(0)\bs\{0\}$$
is a conformal map. In particular, we obtain $f(B_\delta^{\circ}\bs\{z_0\})\subset B_\delta^{\circ}\bs\{z_0\}$. Endowing the set $\mathcal{K}(B_\delta)$ with the Hausdorff distance $d_{\mathcal{K}(B_\delta)}$, we obtain the complete metric space $\left(\mathcal{K}(B_\delta),d_{\mathcal{K}(B_\delta)}\right)$ (see e.g.\,\cite{edgar}, p.66f.). In the following, we want to show that $C_f$ is $A(K)$-universal for comeager many $K\in\mathcal{K}(B_\delta)$. The starting point is the following universality result, which holds for finite sets:
\begin{cor}\label{superattracting_finite}
Let $E\subset B_\delta^{\circ}\bs\{z_0\}$ be a finite set such that all iterates $f^n$ are injective on $E$. Then the set of all functions in $H(\mathbb{C}\bs\{z_0\})$ which are $C(E)$-universal for $C_f$ is a comeager set in $H(\mathbb{C}\bs\{z_0\})$.
\end{cor}
\textbf{Proof:}
We consider the domain $\Omega:=\mathbb{C}\bs\{z_0\}$ and the open set $D:=B_\delta^{\circ}\bs\{z_0\}\subset\Omega$. Due to the injectivity of all iterates $f^n|_E$ and the convergence $f^n|_E\rightarrow z_0\in\partial_\infty\Omega$, we obtain $E\in\mathcal{E}(D,\Omega,f)$ so that the assertion now follows from Corollary \ref{central_universality_theorem_finite}.\hfill$\square$\medskip

In order to extend the universality statement of Corollary \ref{superattracting_finite} to comeager many compact subsets of $B_\delta^{\circ}\bs\{z_0\}$, we denote $A_\delta:=(B_\delta^\circ\bs\{z_0\})\cap(\mathbb{Q}+i\mathbb{Q})$ and we introduce the countable set
$$\mathcal{E}_\delta:=\left\{E\subset A_\delta:E\neq\emptyset,\,E\mbox{ finite, }f^n|_E\mbox{ injective for all }n\in\mathbb{N}\right\}.$$
\begin{lem}\label{dense}
$\mathcal{E}_\delta$ is dense in $\big(\mathcal{K}(B_\delta),d_{\mathcal{K}(B_\delta)}\big)$ and \,$\mathcal{K}(B_\delta^\circ\bs\{z_0\})$ is open and dense in $\big(\mathcal{K}(B_\delta),d_{\mathcal{K}(B_\delta)}\big)$.
\end{lem}
\textbf{Proof:} The proof of the first statement is rather technical -- therefore, we refer to Lemma 3.3.2 and Corollary 3.3.3 in \cite{jung}. In order to prove the second statement, let $K\in\mathcal{K}(B_\delta^\circ\bs\{z_0\})$. Putting $\delta:=\dist(\{z_0\}\cup\partial B_\delta,K)>0$, we obtain
$\big\{L\in\mathcal{K}(B_\delta):d_{\mathcal{K}(B_\delta)}(K,L)<\delta\big\}\subset\mathcal{K}(B_\delta^\circ\bs\{z_0\}).$  
Indeed, let $L\in \mathcal{K}(B_\delta)$ with $d_{\mathcal{K}(B_\delta)}(K,L)<\delta$. Assuming that there exists a point $\wt{z}\in L\cap(\{z_0\}\cup\partial B_\delta)$, we would obtain
$$d_{\mathcal{K}(B_\delta)}(K,L)\geq\max_{z\in L}\dist(z,K)\geq\dist(\wt{z},K)\geq\dist(\{z_0\}\cup\partial B_\delta,K)=\delta,$$
a contradiction. Thus, we have $L\cap(\{z_0\}\cup\partial B_\delta)=\emptyset$ and hence $L\in\mathcal{K}(B_\delta^\circ\bs\{z_0\})$. Due to $\mathcal{K}(B_\delta^\circ\bs\{z_0\})\supset\mathcal{E}_\delta$, the first part implies the denseness of $\mathcal{K}(B_\delta^\circ\bs\{z_0\})$ in $\big(\mathcal{K}(B_\delta),d_{\mathcal{K}(B_\delta)}\big)$.\hfill$\square$\medskip

Now, we can extend the statement of Corollary \ref{superattracting_finite} in the following way:
\begin{theo}\label{superattracting_1}
For each $\delta>0$ with $U_\delta[0]\subset V,$ comeager many functions in $H(\mathbb{C}\bs\{z_0\})$ are $C(K)$-universal for $C_f$ for comeager many $K\in\mathcal{K}(B_\delta)$.
\end{theo}
\textbf{Proof:} The proof will run similarly to the proof of Lemma 2 in \cite{mueller}.
\begin{itemize}
\item[i)]For $g\in H(\mathbb{C}\bs\{z_0\})$, we define
$$\mathcal{K}_g:=\Big\{K\in\mathcal{K}(B^\circ_\delta\bs\{z_0\}):\left\{g\circ f^n|_K:n\in\mathbb{N}\right\}\mbox{ dense in }C(K)\Big\}.$$
Denoting by $\mathcal{P}$ the set of all complex-valued polynomials in two real variables with Gaussian rational coefficients,
the complex Stone-Weierstrass theorem implies that the set $\left\{p|_L:p\in\mathcal{P}\right\}$ is dense in $C(L)$ for each $L\in\mathcal{K}(\mathbb{C})$. Hence, we obtain
$$\mathcal{K}_g=\bigcap_{\substack{j\in\mathbb{N}\\p\in\mathcal{P}}}\bigcup_{n\in\mathbb{N}}\left\{K\in\mathcal{K}(B^\circ_\delta\bs\{z_0\}):\left\|g\circ f^n|_K-p\right\|_K<\frac{1}{j}\right\}.$$
\item[ii)]Let $g\in H(\mathbb{C}\bs\{z_0\})$, $p\in\mathcal{P}$ and $n\in\mathbb{N}$ be fixed. We show that the function
$$\psi:\mathcal{K}(B_\delta^\circ\bs\{z_0\})\rightarrow[0,\infty),\,\,\,\psi(K):=\left\|g\circ f^n|_K-p\right\|_K$$
is continuous. In order to do so, we denote $q:=g\circ f^n|_{B_\delta^\circ\bs\{z_0\}}-p$. 
Now, let $K\in\mathcal{K}(B_\delta^\circ\bs\{z_0\})$ and $\varepsilon>0$. We choose some $\delta_1>0$ with $U_{\delta_1}[K]:=\{z\in\mathbb{C}:\dist(z,K)\leq\delta_1\}\subset B_\delta^\circ\bs\{z_0\}$. As $q$ is uniformly continuous on $U_{\delta_1}[K]$, there exists a $\delta_2>0$ such that we have
$|q(z_1)-q(z_2)|<\varepsilon$
for all points $z_1,z_2\in U_{\delta_1}[K]$ with $|z_1-z_2|<\delta_2$. We put $\delta:=\min(\delta_1,\delta_2)$ and we consider a set $L\in\mathcal{K}(B_\delta^\circ\bs\{z_0\})$ with $d_{\mathcal{K}(B_\delta)}(L,K)<\delta$. Then, for $z\in L$, there exists a point $w\in K$ with $|z-w|<\delta\leq\delta_2$, and we have $\dist(z,K)\leq |z-w|<\delta$ which yields $z\in U_{\delta_1}[K]$. Therefore, we obtain
$|q(z)|\leq|q(z)-q(w)|+|q(w)|<\varepsilon+\|q\|_K.$
As this inequality holds for all $z\in L$, we obtain $\|q\|_L\leq\varepsilon+\|q\|_K$ and hence $\|q\|_L-\|q\|_K\leq\varepsilon$. Analogously, it follows that $\|q\|_K-\|q\|_L\leq\varepsilon$, which finally yields
$$|\psi(L)-\psi(K)|=|\|q\|_L-\|q\|_K|\leq\varepsilon.$$
Thus, each set
$\big\{K\in\mathcal{K}(B_\delta^\circ\bs\{z_0\}):\left\|g\circ f^n|_K-p\right\|_K<1\slash j\big\}$
is open in $\mathcal{K}(B_\delta^\circ\bs\{z_0\})$. Due to Lemma \ref{dense}, these sets are also open in $\mathcal{K}(B_\delta)$. Hence, part i) yields that $\mathcal{K}_g$ is a $G_\delta$-set in $\mathcal{K}(B_\delta)$.
\item[iii)]For each $E\in\mathcal{E}_\delta$, Corollary \ref{superattracting_finite} implies that the set
$$\mathcal{G}_E:=\big\{g\in H(\mathbb{C}\bs\{z_0\}):g\,\,C(E)\mbox{-universal for } C_f\big\}$$ 
is comeager in $H(\mathbb{C}\bs\{z_0\})$. As $\mathcal{E}_\delta$ is countable, $\mathcal{G}:=\bigcap_{E\in\mathcal{E}_\delta}\mathcal{G}_E$
is also comeager in $H(\mathbb{C}\bs\{z_0\})$. Moreover, as each function $g\in\mathcal{G}$ is $C(E)$-universal for $C_f$ for all sets $E\in\mathcal{E}_\delta$, we have $\mathcal{E}_\delta\subset\mathcal{K}_g$ for all $g\in\mathcal{G}$. As $\mathcal{E}_\delta$ is dense in $\mathcal{K}(B_\delta)$ due to Lemma \ref{dense}, it follows that $\mathcal{K}_g$ is dense in $\mathcal{K}(B_\delta)$ for all $g\in\mathcal{G}$. Therefore, part ii) yields that for comeager many $g\in H(\mathbb{C}\bs\{z_0\})$ the set $\mathcal{K}_g$ is a dense $G_\delta$-set in $\mathcal{K}(B_\delta)$ and hence, in particular, comeager in $\mathcal{K}(B_\delta)$. Thus, comeager many functions $g\in H(\mathbb{C}\bs\{z_0\})$ have the property that for comeager many $K\in\mathcal{K}(B_\delta)$ the set $\left\{g\circ f^n|_K:n\in\mathbb{N}\right\}$ is dense in $C(K)$.\hfill$\square$
\end{itemize}\medskip

It is well-known that comeager many sets in $\mathcal{K}(B_\delta)$ are Cantor sets, i.e.\,\,perfect and totally disconnected (cf.\,\cite{bmm}, Remark 2 on p.\,236). In particular, we have $K^\circ=\emptyset$ and hence $A(K)=C(K)$ for comeager many $K\in\mathcal{K}(B_\delta)$. For this reason, we may replace $C(K)$-universality by $A(K)$-universality in the statement of Theorem \ref{superattracting_1}. However, as the considerations before Corollary \ref{superattracting_finite} indicate, for any set $K\in\mathcal{K}(B_\delta^\circ\bs\{z_0\})$ with $K^\circ\neq\emptyset$, there does not exist a function in $H(U\bs\{z_0\})$ which is $A(K)$-universal for $C_f$ (cf.\,\cite{jung}, Remark 3.3.6\,iii)). 
In particular, the universality statement of Theorem \ref{superattracting_1} is much weaker than the universality result which holds near attracting fixed points of $f$ (Theorem \ref{attracting_case}), where we actually obtained $H(W)$-universality of $C_f$ for suitable open sets $W$.\medskip

We want to conclude this subsection with the following consideration. For $\delta>0$ with $U_\delta[0]\subset V$, Theorem \ref{superattracting_1} states that comeager many functions in $H(\mathbb{C}\bs\{z_0\})$ are $C(K)$-universal for $C_f$ for comeager many $K\in\mathcal{K}(B_\delta)$. Here, the second ``comeager many''-expression depends on the first one, meaning that the comeager subset of $\mathcal{K}(B_\delta)$ depends on the choice of a function from the comeager subset of $H(\mathbb{C}\bs\{z_0\})$ (cf.\,\,the proof of Theorem \ref{superattracting_1}). In the following, we want to show that this dependence can be interchanged. For two sets $X$ and $Y$ and $x_0\in X$, $y_0\in Y$ as well as $A\subset X\times Y$, we denote
$$A(x_0,\cdot):=\{y\in Y:(x_0,y)\in A\}\quad\mbox{and}\quad A(\cdot,y_0):=\{x\in X:(x,y_0)\in A\}.$$
\begin{lem}\label{quasi_quasi}  
Let $X$ and $Y$ be Baire spaces and let $Y$ be second-countable. Moreover, let $A\subset X\times Y$ such that $A(x,\cdot)$ is a $G_\delta$-set in $Y$ for comeager many $x\in X$ and such that $A(\cdot,y)$ is comeager in $X$ for comeager many $y\in Y$. Then $A(x,\cdot)$ is comeager in $Y$ for comeager many $x\in X$.
\end{lem}
\textbf{Proof:} According to the assumption, there exists a comeager subset $S$ of $Y$ such that $A(\cdot,y)$ is comeager in $X$ for all $y\in S$. As $Y$ is a Baire space, we obtain that $S$ is dense in $Y$. The second-countability of $Y$ implies that $S$ is also second-countable and hence, in particular, separable. Thus, there exists a countable subset $R$ of $S$ which is dense in $S$. Then $R$ is also dense in $Y$ and the set
$X_0:=\bigcap_{y\in R}A(\cdot,y)$
is comeager in $X$. For $x_0\in X_0$, we have $(x_0,y)\in A$ and hence $y\in A(x_0,\cdot)$ for all $y\in R$. Thus, we obtain $A(x_0,\cdot)\supset R$ so that the denseness of $R$ in $Y$ implies that $A(x_0,\cdot)$ is dense in $Y$. Hence, $A(x,\cdot)$ is dense in $Y$ for comeager many $x\in X$. The assertion now follows from the assumption that $A(x,\cdot)$ is also a $G_\delta$-set in $Y$ for comeager many $x\in X$.\hfill$\square$\medskip

There exists a similar version of Lemma \ref{quasi_quasi} which follows from the Kuratowski-Ulam theorem (see e.g.\,\cite{kahane}, p.\,144). This version states that for a Baire space $X$, a second-countable Baire space $Y$ and a set $A\subset X\times Y$ the following is true: If $A$ is a $G_\delta$-set in $X\times Y$ such that $A(\cdot,y)$ is comeager in $X$ for all $y\in Y$, then $A(x,\cdot)$ is comeager in $Y$ for comeager many $x\in X$. As the property of $A$ being a $G_\delta$-set in $X\times Y$ implies that $A(x,\cdot)$ is a $G_\delta$ in $Y$ for all $x\in X$, Lemma \ref{quasi_quasi} yields the same result as the above-mentioned version -- but under weaker assumptions. In order to prove the following theorem, we will see that we do indeed need the statement of Lemma \ref{quasi_quasi}.

\begin{theo}\label{superattracting_2}
For each $\delta>0$ with $U_\delta[0]\subset V,$ comeager many $K\in\mathcal{K}(B_\delta)$ have the property that comeager many functions in $H(\mathbb{C}\bs\{z_0\})$ are $C(K)$-universal for $C_f$.
\end{theo}
\textbf{Proof:} We denote $X:=\mathcal{K}(B_\delta)$, $Y:=H(\mathbb{C}\bs\{z_0\})$ and
$$A:=\Big\{(K,g)\in\mathcal{K}(B_\delta^\circ\bs\{z_0\})\times H(\mathbb{C}\bs\{z_0\}):\left\{g\circ f^n|_K:n\in\mathbb{N}\right\}\mbox{ dense in }C(K)\Big\}.$$
Then $X$ and $Y$ are Baire spaces, $Y$ is second-countable and $A$ is a subset of \mbox{$X\times Y$}. Denoting by $\mathcal{P}$ the set of all complex-valued polynomials in two real variables with Gaussian rational coefficients, the complex Stone-Weierstrass theorem implies, \mbox{analogously} to the proof of Theorem \ref{superattracting_1}, that we have
\begin{eqnarray*}
A(K,\cdot)&=&\Big\{g\in H(\mathbb{C}\bs\{z_0\}):\left\{g\circ f^n|_K:n\in\mathbb{N}\right\}\mbox{ dense in }C(K)\Big\}\\[2pt]
&=&\bigcap_{\substack{j\in\mathbb{N}\\p\in\mathcal{P}}}\bigcup_{n\in\mathbb{N}}\left\{g\in H(\mathbb{C}\bs\{z_0\}):\|g\circ f^n|_K-p\|_K<\frac{1}{j}\right\}
\end{eqnarray*}
for all $K\in\mathcal{K}(B_\delta^\circ\bs\{z_0\})$. For each $p\in\mathcal{P},n\in\mathbb{N}$ and $K\in\mathcal{K}(B_\delta^\circ\bs\{z_0\})$, the continuity of the composition operator $C_{f^n,K}:H(\mathbb{C}\bs\{z_0\})\rightarrow A(K),\,C_{f^n,K}(g):=g\circ f^n|_K$, yields that the function
$$H(\mathbb{C}\bs\{z_0\})\ni g\mapsto\left\|g\circ f^n|_K-p\right\|_K$$
is continuous. Hence, each set
$\big\{g\in H(\mathbb{C}\bs\{z_0\}):\|g\circ f^n|_K-p\|_K<1\slash j\big\}$
is open in $H(\mathbb{C}\bs\{z_0\})$ so that we obtain that $A(K,\cdot)$ is a $G_\delta$-set in $Y$ for all $K\in\mathcal{K}(B_\delta^\circ\bs\{z_0\})$. As $\mathcal{K}(B_\delta^\circ\bs\{z_0\})$ is a comeager subset of $\mathcal{K}(B_\delta)$ due to Lemma \ref{dense}, it follows that $A(K,\cdot)$ is a $G_\delta$-set in $Y$ for comeager many $K\in X$. According to Theorem \ref{superattracting_1}, we have that
$$A(\cdot,g)=\Big\{K\in\mathcal{K}(B_\delta^\circ\bs\{z_0\}):\left\{g\circ f^n|_K:n\in\mathbb{N}\right\}\mbox{ dense in }C(K)\Big\}$$
is comeager in $X$ for comeager many $g\in Y$ (see also the proof of Theorem \ref{superattracting_1}, where we have $A(\cdot,g)=\mathcal{K}_g$). Therefore, Lemma \ref{quasi_quasi} implies that for comeager many $K\in X$ the set $A(K,\cdot)$ is comeager in $Y$.\hfill$\square$

\section{Applications: Global Theory}    
\label{sec:Applications:Global Theory} 

For a complex-valued function $f$ which is holomorphic on an open neighbourhood of a (super-)attracting or a neutral fixed point $z_0\in\mathbb{C}$ of $f$, Section \ref{sec:Applications: Local Theory} provides several universality statements for $C_f$. If, in addition, the symbol $f$ is not only locally defined but even an entire function (other than a polynomial of degree $0$ or $1$), all these universality phenomenons occur on ``small'' parts of the Fatou set of $f$. (Note that (super-)attracting fixed points of $f$ as well as attracting petals for $f$ at neutral fixed points are always contained in $F_f$.) More precisely, this means that we were able to find open sets $U\subset F_f$ or compact sets $K\subset F_f$ located near $z_0$ with the property that $C_f$ is $H(U)$- or $A(K)$-universal, respectively. In case of $H(U)$-universality of $C_f$, there exist (comeager many) functions $g\in H(\mathbb{C}\bs\{z_0\})$ for which we have
$$\ol{\left\{g\circ f^n|_U:n\in\mathbb{N}\right\}}^{H(U)}=H(U).$$
In this section, we will study the long-time behaviour of sequences of compositions of the form $(g\circ f^n)$ on ``large'' parts of the Fatou set of an entire symbol $f$ (other than a polynomial of degree $0$ or $1$), e.g.\,\,on whole components of $F_f$. Our aim is to determine the ``richness'' of the closure of sets of the form $\left\{g\circ f^n|_G:n\in\mathbb{N}\right\}$ for ``many'' functions $g\in H(\Omega)$, where now $G$ should be a ``large'' open invariant subset of $F_f$ and $\Omega$ is a domain in $\mathbb{C}$ which contains $G$. In contrast to the local situations which we have considered in the previous section, we now cannot expect that there exists a function $g\in H(\Omega)$ such that the closure of the set of compositions $\left\{g\circ f^n|_G:n\in\mathbb{N}\right\}$ equals the whole space $H(G)$, i.e.\,\,we cannot expect that $C_f$ is $H(G)$-universal. The reason for this is that whenever $f$ is not injective on $G$ -- which is likely to be the case for ``large'' subsets $G$ of $F_f$ -- then there exists the following ``natural restriction'' for each function \label{universal_implies_injective}
$$h\in\ol{\left\{g\circ f^n|_G:n\in\mathbb{N}\right\}}^{H(G)}\,\big\bs\,\left\{g\circ f^n|_G:n\in\mathbb{N}\right\}$$
(cf.\,\,the considerations after Theorem \ref{theorem_grosse_erdmann_mortini}): Choosing a strictly increasing sequence $(n_k)$ in $\mathbb{N}$ such that $(g\circ f^{n_k})$ converges to $h$ locally uniformly on $G$, we obtain for any $N\in\mathbb{N}$ and any points $z,w\in G$ with $z\neq w$ and $f^N(z)=f^N(w)$ that
$$h(z)=\lim_{k\rightarrow\infty}g(f^{n_k}(z))=\lim_{k\rightarrow\infty}g(f^{n_k}(w))=h(w).$$
This means that we can approximate only those functions in $H(G)$ via subsequences of $(g\circ f^n)$ which have the property of assuming the same value at all points in $G$ which eventually coincide under the iteration of $f$. 
Therefore, the composition operator $C_f$ cannot be $H(G)$-universal. In this situation, the challenge now consists in determining 
the closure of the set of compositions $\left\{g\circ f^n|_G:n\in\mathbb{N}\right\}$ in $H(G)$. In order to do so, it is useful to consider the set $\omega(G,g,f)$, which shall denote the set of all functions $h:G\rightarrow\mathbb{C}$ for which there exists a strictly increasing sequence $(n_k)$ in $\mathbb{N}$ such that $(g\circ f^{n_k})$ converges to $h$ locally uniformly on $G$. Using this terminology, we have
$$\ol{\left\{g\circ f^n|_G:n\in\mathbb{N}\right\}}^{H(G)}=\,\,\omega(G,g,f)\,\cup\,\left\{g\circ f^n|_G:n\in\mathbb{N}\right\}.$$

\subsection{Schr\"oder Domains}
\label{subsec:Schröder Domains}

Let $z_0\in\mathbb{C}$ be an attracting fixed point of $f$ and let
$A(z_0):=\left\{z:f^n(z)\rightarrow z_0\right\}$
be the \textit{basin of attraction} of $z_0$ under $f$, which is a union of components of $F_f$ (see e.g.\,\cite{carleson-gamelin}, Theorem III.2.1). In particular, the Schr\"oder domain of $f$ containing $z_0$ (i.e.\,\,the component of $F_f$ which contains $z_0$) is a subset of $A(z_0)$ (cf.\,\,Theorem \ref{classifcation_theorem}). Moreover, we consider the \textit{backward orbit} of $z_0$ under $f$, which is given by 
$O^-(z_0):=\bigcup_{n\in\mathbb{N}}\left\{z:f^n(z)=z_0\right\}.$
By definition, the set $A(z_0)\bs O^-(z_0)$ is completely invariant under $f$. In particular, $f$ is not injective on \,$A(z_0)\bs O^-(z_0)$.\medskip

As already described in the introduction of this chapter, it is now our aim to characterize sets of the form $\omega(G,g,f)$ for suitable ``large'' open subsets $G$ of $A(z_0)\bs O^-(z_0)$ and ``many'' functions $g\in H(\mathbb{C}\bs\{z_0\})$. Due to the above considerations, each function $h\in\omega(G,g,f)$ must have the property of assuming the same value at all points $z,w\in G$ for which we have $f^N(z)=f^N(w)$ for some $N\in\mathbb{N}$. Thus, in order to determine the set $\omega(G,g,f)$, we have to look for holomorphic functions on $G$ which fulfil this property. Clearly, this holds for all constant functions -- but the crucial question we are facing here is whether there exist any non-constant functions for which this is true. As we do not have any information concerning this fact at first view, it is natural to remind ourselves that $f$ is locally conjugated near the attracting fixed point $z_0$ to the linear function $w\mapsto\lambda w$ with $\lambda:=f'(z_0)$ (where $0<|\lambda|<1$).  
This means that there exist open neighbourhoods $U$ of $z_0$ and $V$ of $0$ with $f(U)\subset U$ and $\lambda V\subset V$ as well as a conformal map $\varphi:U\rightarrow V$ such that the equation     
$\varphi\circ f^n=\lambda^n\cdot\varphi$
holds on $U$ for all integers $n\in\mathbb{N}$ (cf.\,\,Subsection \ref{subsec:Attracting Fixed Points}). According to the considerations before Theorem \ref{attracting_case}, we have $U\subset A(z_0)$. It is well-known that the conjugation map $\varphi$ can be extended holomorphically to the entire basin of attraction of $z_0$ and that the extended $\varphi$ fulfils the same functional equation as $\varphi$ (see e.g.\,\cite{carleson-gamelin}, p.\,32), i.e.\,\,there exists a holomorphic function $\Phi:A(z_0)\rightarrow\mathbb{C}$ with $\Phi|_U=\varphi$ such that the equation $\Phi\circ f^n=\lambda^n\cdot\Phi$ holds on $A(z_0)$ for all $n\in\mathbb{N}$. As the image $\Phi(A(z_0))$ is a dense subset of the complex plane (cf.\,\cite{jung}, Remark 4.1.1 and p.\,66), the function $\Phi$ quasiconjugates $f|_{A(z_0)}$ to $w\mapsto\lambda w$ on $\mathbb{C}$. For $z,w\in A(z_0)$ and $N\in\mathbb{N}$ with $f^N(z)=f^N(w)$, we have
$$\Phi(z)=\frac{\Phi(f^N(z))}{\lambda^N}=\frac{\Phi(f^N(w))}{\lambda^N}=\Phi(w),$$
i.e.\,\,the non-constant holomorphic function $\Phi$ assumes the same value at all points in $A(z_0)$ which eventually coincide under the iteration of $f$. 
Knowing that $\Phi$ fulfils the property of assuming the same value at all points in $A(z_0)$ which eventually coincide under the iteration of $f$, we obtain that this is already true for a whole class of functions, namely for all compositions $\psi\circ\Phi$, where $\psi$ is defined on $\Phi(A(z_0))$. Therefore, we introduce the following notation:
\begin{defi}\label{defi_h_phi}
For $G\subset A(z_0)$ open, we put
$$H_{\Phi}(G):=\left\{\psi\circ\Phi|_G:\psi\in H(\Phi(G))\right\}.$$
\end{defi}

We remark that the set $H_{\Phi}(G)$ is a subspace of $H(G)$  and that it does not depend on the choice of the function $\Phi$, i.e.\,\,if there exist open neighbourhoods $U'$ of $z_0$ and $V'$ of $0$ as well as a conformal map $\wt{\varphi}:U'\rightarrow V'$ with $\wt{\varphi}\circ f=\lambda\cdot\wt{\varphi}$ on $U'$, then there exists some $c\in\mathbb{C}\bs\{0\}$ with $\varphi=c\cdot\wt{\varphi}$ on $U\cap U'$ (cf.\,\cite{carleson-gamelin}, p.\,28f.). Hence, the holomorphic extension $\wt{\Phi}$ of $\wt{\varphi}$ to $A(z_0)$ which fulfils $\wt{\Phi}\circ f=\lambda\cdot\wt{\Phi}$ on $A(z_0)$ equals $c\cdot\wt{\Phi}$ so that we obtain $H_{\wt{\Phi}}(G)=H_{\Phi}(G)$.\medskip

For each function $h\in H_{\Phi}(G)$, we have $h(z)=h(w)$ for all points $z,w\in G$ with $f^N(z)=f^N(w)$ for some $N\in\mathbb{N}$. Thus, for \,$G\subset A(z_0)\bs O^-(z_0)$ open and $g\in H(\mathbb{C}\bs\{z_0\})$, each function in $H_{\Phi}(G)$ is a possible candidate for being contained in the set $\omega(G,g,f)$.
In the following, it will be our aim to prove the existence of a ``large'' open subset $G_0$ of $A(z_0)\bs O^-(z_0)$ such that we actually have\,
$\omega(G_0,g,f)=H_{\Phi}(G_0)$
and hence, in particular,
$$\ol{\big\{g\circ f^n|_{G_0}:n\in\mathbb{N}\big\}}^{H(G_0)}=\,\,H_{\Phi}(G_0)\,\cup\,\big\{g\circ f^n|_{G_0}:n\in\mathbb{N}\big\}$$
for ``many'' functions $g\in H(\mathbb{C}\bs\{z_0\})$. In order to show that such a ``large'' set $G_0$ exists, we introduce the following definition: 
\begin{defi}\label{locally_universal}
A function $g\in H(A(z_0)\bs\{z_0\})$ is called \textit{locally universal} for $C_f$  
if there exist open neighbourhoods $U$ of $z_0$ and $V$ of $0$ as well as a conformal map $\varphi:U\rightarrow V$ with $\varphi\circ f=\lambda\cdot\varphi$ on $U$ such that
the set $\{g\circ f^n|_W:n\in\mathbb{N}\}$ is dense in $H(W)$ for all $W\in\mathcal{U}_0(U\bs\{z_0\})$.
\end{defi}
As Picard's little theorem implies either $f^n(\mathbb{C})=\mathbb{C}$ or the existence of a point $a_f\in\mathbb{C}$ with $f^n(\mathbb{C})=\mathbb{C}\bs\{a_f\}$ for all $n\in\mathbb{N}$, we may assume without loss of generality that we always have $a_f\notin U$ 
whenever we consider a locally universal function $g$ for $C_f$. 
Theorem \ref{attracting_case} yields that the set
$$\mathcal{G}:=\left\{g\in H(\mathbb{C}\bs\{z_0\}):g|_{A(z_0)\bs\{z_0\}}\mbox{ locally universal for }C_f\right\}$$
is a comeager set in $H(\mathbb{C}\bs\{z_0\})$. 
The following lemma generalizes the statement of Theorem \ref{attracting_case}:
\begin{lem}\label{attracting_case_generalized}
Let $g$ be locally universal for $C_f$, $N\in\mathbb{N}_0$ and $W\in\mathcal{U}_0(U\bs\{z_0\})$. Then we have
$\omega\big(f^{-N}(W),g,f\big)\supset H_\Phi\big(f^{-N}(W)\big).$  
\end{lem}
\textbf{Proof:} We first observe that we have $f^N(f^{-N}(W))=W$. In fact, the inclusion $f^N(f^{-N}(W))\subset W$ always holds, and the reverse inclusion is true whenever the set $W$ is contained in the range of $f^N$. As the possibly existing Picard exceptional value of $f$ is not contained in $U\supset W$, this is indeed the case. Hence, as $\Phi$ quasi-conjugates $f|_{A(z_0)}$ to $w\mapsto\lambda w$ on $\mathbb{C}$, we have
$$\Phi\big(f^{-N}(W)\big)=\frac{1}{\lambda^N}\cdot\Phi\big(f^N(f^{-N}(W))\big)=\frac{1}{\lambda^N}\cdot\Phi(W).$$
Now, let $h\in H_\Phi\big(f^{-N}(W)\big)$. Then there exists some $\psi\in H\big(\Phi(f^{-N}(W))\big)$ with 
$$h=\psi\circ\Phi|_{f^{-N}(W)}=\psi\circ\frac{1}{\lambda^N}\cdot\Phi\circ f^N|_{f^{-N}(W)}=\wt{h}\circ f^N|_{f^{-N}(W)},$$
where $\wt{h}:=\psi\circ\frac{1}{\lambda^N}\cdot\Phi|_W\in H(W).$ 
For $K\subset f^{-N}(W)$ compact, $f^N(K)$ is a compact subset of $W$ so that the denseness of $\left\{g\circ f^n|_W:n\in\mathbb{N}\right\}$ in $H(W)$ yields an integer $m_1\in\mathbb{N}$ with
$$1>\left\|g\circ f^{m_1}|_W-\wt{h}\right\|_{f^N(K)}=\|g\circ f^{m_1}\circ f^N-\wt{h}\circ f^N\|_K=\|g\circ f^{m_1+N}-h\|_K.$$
As the space $H(W)$ has no isolated points, each set $\left\{g\circ f^n|_W:n\geq m\right\}$ is dense in $H(W)$. Thus, we inductively find a strictly increasing sequence $(m_j)$ in $\mathbb{N}$ with 
$\|g\circ f^{m_j+N}-h\|_K<1\slash j$ 
for all $j\in\mathbb{N}$. Hence, we have uniform convergence $g\circ f^{m_j+N}\rightarrow h$ on $K$. But this already implies the existence of a strictly increasing sequence $(n_k)$ in $\mathbb{N}$ with $g\circ f^{n_k}\rightarrow h$ locally uniformly on $f^{-N}(W)$ (see e.g.\,\cite{jung}, Lemma A.4), i.e.\,\,we have $h\in\omega\big(f^{-N}(W),g,f\big)$.\hfill$\square$\medskip

For $g\in\mathcal{G}$, $W\in\mathcal{U}_0(U\bs\{z_0\})$ and $N=0$, Lemma \ref{attracting_case_generalized} states $\omega(W,g,f)\supset H_\Phi(W)$. As each function $h\in H(W)$ can be written as
$$h=h\circ\varphi^{-1}|_{\varphi(W)}\circ\varphi|_W=\left(h\circ\varphi^{-1}|_{\varphi(W)}\right)\circ\Phi|_W\in H_\Phi(W),$$
we actually have $H_\Phi(W)=H(W)$ and hence $\omega(W,g,f)\supset H(W)$, which is equivalent to the fact that the set $\left\{g\circ f^n|_W:n\in\mathbb{N}\right\}$ is dense in $H(W)$. Therefore, the statement of Lemma \ref{attracting_case_generalized} in case of $g\in\mathcal{G}$ and $N=0$ is just the statement of Theorem \ref{attracting_case}.

\begin{cor}\label{attracting_case_first_inlcusion}
Let $g$ be locally universal for $C_f$ and let $W\in\mathcal{U}_0(U\bs\{z_0\})$ with $f(W)\subset W$. Then 
$D:=\bigcup_{n\in\mathbb{N}_0}f^{-n}(W)$
is an open subset of $A(z_0)\bs O^-(z_0)$ and we have $\omega(D,g,f)\supset H_\Phi(D)$. 
\end{cor}
\textbf{Proof:} Because of $W\subset U\subset A(z_0)$, the backward invariance of $A(z_0)$ under $f$ yields $D\subset A(z_0)$. Assuming that there exists a point $z\in D\cap O^-(z_0)$, there would exist integers $n\in\mathbb{N}_0$, $k\in\mathbb{N}$ with $f^n(z)\in W$ and $f^k(z)=z_0$. In case of $k\leq n$, we would obtain $z_0=f^n(z)\in W\subset U\bs\{z_0\}$, a contradiction. For $k>n$, it would follow that $z_0=f^k(z)=f^{k-n}(f^n(z))\in f^{k-n}(W)\subset W\subset U\bs\{z_0\}$, the same contradiction. Thus, we have $D\subset A(z_0)\bs O^-(z_0)$.
Now, let $h\in H_\Phi(D)$. Then there exists some $\psi\in H(\Phi(D))$ with $h=\psi\circ\Phi|_{D}$. For $K\subset D$ compact, there exists a finite set $E\subset\mathbb{N}_0$ with $K\subset\bigcup_{n\in E}f^{-n}(W)$. Putting $N:=\max E$, the invariance of $W$ under $f$ yields
$$f^N(K)\subset\bigcup_{n\in E}f^N(f^{-n}(W))\subset\bigcup_{n\in E}f^{N-n}(W)\subset W.$$
Hence, $K$ is a compact subset of $f^{-N}(W)$. According to Lemma \ref{attracting_case_generalized}, we have
$$h|_{f^{-N}(W)}=\psi|_{\Phi\left(f^{-N}(W)\right)}\circ\Phi|_{f^ {-N}(W)}\in H_\Phi\big(f^{-N}(W)\big)\subset\omega\big(f^{-N}(W),g,f\big).$$
Thus, there exists a strictly increasing sequence $(m_j)$ in $\mathbb{N}$ such that $(g\circ f^{m_j})$ converges to $h$ locally uniformly on $f^{-N}(W)$. In particular, we obtain uniform convergence $g\circ f^{m_j}\rightarrow h$ on $K$. As above, this implies the existence of a strictly increasing sequence $(n_k)$ in $\mathbb{N}$ with $g\circ f^{n_k}\rightarrow h$ locally uniformly on $D$, i.e.\,\,we have $h\in\omega(D,g,f)$.\hfill$\square$\medskip

\begin{figure}[h]
\centering
\includegraphics[width=12.323cm,height=6.107cm]{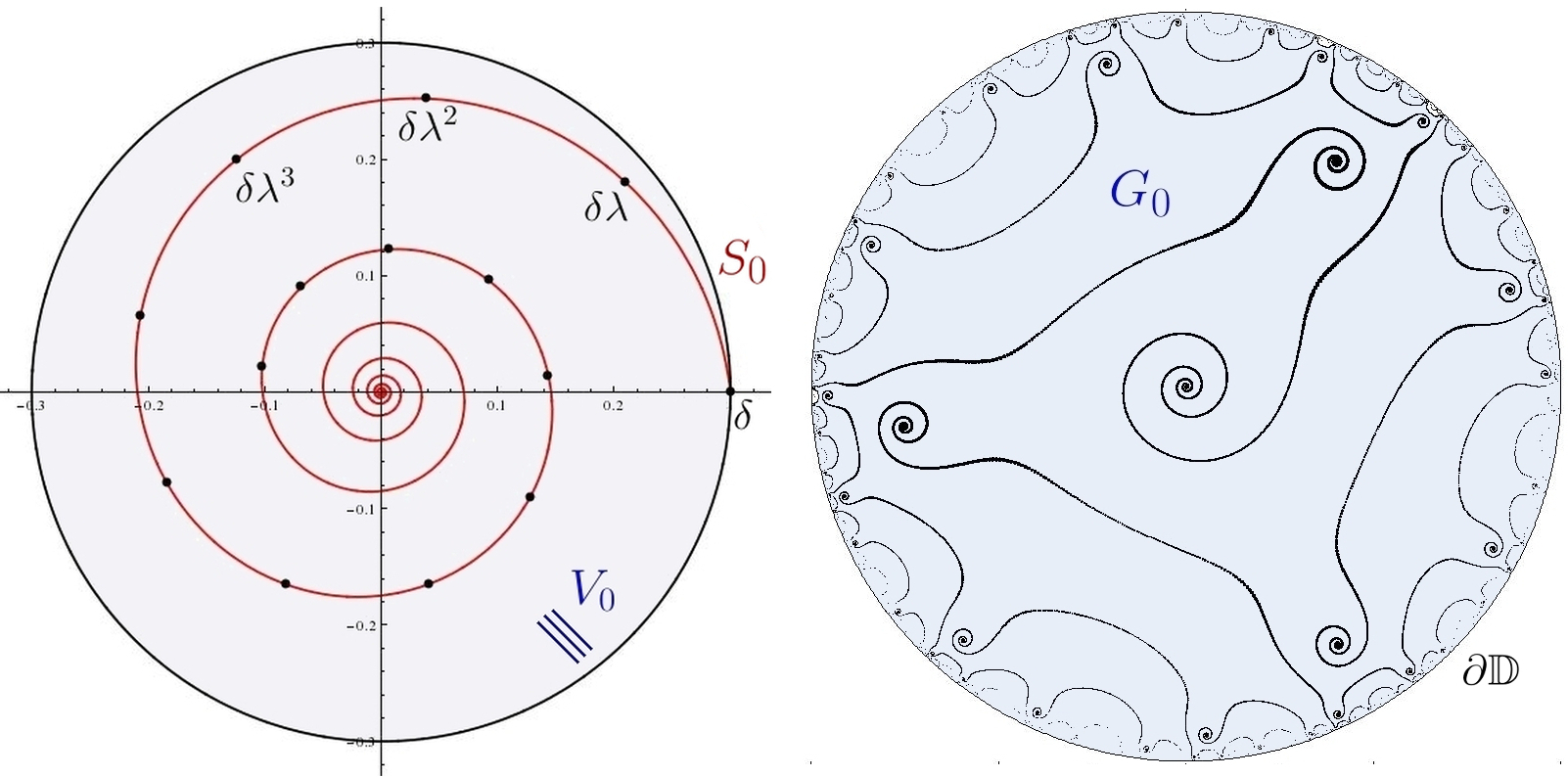}
\caption{}
\label{spiral_and_G_0}
\end{figure}

In the following lemma, we prove the existence of a set which fulfils the assumptions of Corollary \ref{attracting_case_first_inlcusion}, i.e.\,\,we construct an open subset $U_0$ of $U\bs\{z_0\}$ which has no holes and which is invariant under $f$. 
\begin{lem}\label{existence_lemma}\quad
\begin{itemize}
\itemsep0pt
\item[i)]Let $W\subset U$ be dense in $U$. Then the set $\bigcup_{n\in\mathbb{N}_0}f^{-n}(W)$ is dense in $A(z_0)$.
\item[ii)]There exists a set $U_0\in\mathcal{U}_0(U\bs\{z_0\})$ with $f(U_0)\subset U_0$ and such that $U_0$ is dense in $U$. 
\end{itemize}
\end{lem}
\textbf{Proof:}\quad
\begin{itemize}
\itemsep0pt
\item[i)]Let $z\in A(z_0)$ and $\varepsilon>0$. Choosing $0<r\leq\varepsilon$ with $U_r[z]\subset A(z_0)$, the Vitali-Porter theorem yields that we have uniform convergence $f^n\rightarrow z_0$ on $U_r[z]$.
Hence, there exists an $N\in\mathbb{N}_0$ with $f^N(U_r[z])\subset U$. As $f^N(U_r(z))$ is a non-empty open subset of $U$, the denseness of $W$ in $U$ implies the existence of a point $w\in W\cap f^N(U_r(z))$ so that there exists a point $\wt{z}\in U_r(z)$ with $f^N(\wt{z})=w\in W$. Thus, we have $|\wt{z}-z|<r\leq\varepsilon$ and $\wt{z}\in f^{-N}(W)\subset\bigcup_{n\in\mathbb{N}_0}f^{-n}(W)$.
\item[ii)]Assuming without loss of generality that $V=U_\delta(0)$ for some $\delta>0$, we first show that there exists a set $V_0\in\mathcal{U}_0(V\bs\{0\})$ with $\lambda V_0\subset V_0$ and such that $V_0$ is dense in $V$. 
In order to do so, we consider the connected set
$S_0:=\left\{\delta\lambda^t:t\in(0,\infty)\right\}\cup\{0\},$
which is the union of the trace of a logarithmic spiral and its ``endpoint'' $0$. 
Putting $V_0:=V\bs S_0,$
we obtain that $V_0$ is an open subset of $V\bs\{0\}$ which is dense in $V$ and has no holes (cf.\,\,the left-hand side of Figure \ref{spiral_and_G_0}).
In order to show that the inclusion $\lambda V_0\subset V_0$ holds, now let $z\in V_0$. Then we have $z\in U_\delta(0)$ and $z\notin S_0$, and due to $0<|\lambda|<1$ it follows that $\lambda z\in V\bs\{0\}$. Assuming that there exists some $t\in(0,\infty)$ with $\lambda z=\delta\lambda^t$, we would obtain $z=\delta\lambda^{t-1}$. On the one hand, this is impossible for $t\leq1$ because in this case it would follow that $|z|=\delta|\lambda|^{t-1}\geq\delta$, a contradiction. On the other hand, $z=\delta\lambda^{t-1}$ is also not possible for $t>1$ because we have $z\notin S_0$. Therefore, we obtain that $\lambda z\in V\bs S_0=V_0$. Hence, there exists a dense subset $V_0$ of $V$ with $V_0\in\mathcal{U}_0(V\bs\{0\})$ and $\lambda V_0\subset V_0$. 
We now consider the set
$U_0:=\varphi^{-1}(V_0),$
which is an open subset of $\varphi^{-1}(V\bs\{0\})=U\bs\{z_0\}$. The denseness of $V_0$ in $V$ and the continuity of $\varphi^{-1}:V\rightarrow U$ imply that $U_0=\varphi^{-1}(V_0)$ is dense in $\varphi^{-1}(V)=U$. Moreover, as $V_0$ does not have any holes, the invariance of the number of holes implies that $U_0$ has no holes as well. Finally, as $\varphi$ conjugates $f|_U$ to $w\mapsto\lambda w$, we obtain from $\lambda V_0\subset V_0$ that
$$f(U_0)=f\big(\varphi^{-1}(V_0)\big)=\varphi^{-1}(\lambda V_0)\subset\varphi^{-1}(V_0)=U_0,$$
which completes the proof.\hfill$\square$
\end{itemize}
\vspace{5pt}

Lemma \ref{existence_lemma} and Corollary \ref{attracting_case_first_inlcusion} yield that $G_0:=\bigcup_{n\in\mathbb{N}_0}f^{-n}(U_0)$ is an open dense subset of $A(z_0)\bs O^-(z_0)$ for which we have $\omega(G_0,g,f)\supset H_\Phi(G_0)$ for each locally universal function $g$ for $C_f$. The right-hand side of Figure \ref{spiral_and_G_0} displays a MATLAB plot of the set $G_0$ in case that $f$ equals the finite Blaschke product
$B:\mathbb{C}_\infty\rightarrow\mathbb{C}_\infty,\,B(z):=z(z-\alpha)\slash(1-\ol{\alpha}z).$
In this situation, we consider the attracting fixed point at the origin (with multiplier $B'(0)=-\alpha$) and its basin of attraction $A(0)=\mathbb{D}$. 
Here, the exceptional set $\mathbb{D}\bs G_0$ looks like a union of pairwise disjoint traces of injective rectifiable paths. In general, it is well-known that the Hausdorff dimension of the trace of an injective rectifiable path defined on a non-singleton compact interval equals one (see e.g.\,\cite{edgar}, Theorem 6.2.7). Indeed, one can show that the Hausdorff dimension of $A(z_0)\bs G_0$ is always one (cf.\,\cite{jung}, Lemma 4.1.9 and Lemma B.7).\medskip

In the following, we will show that the reverse inclusion $\omega(D,g,f)\subset H_\Phi(D)$ holds for each open subset $D$ of $A(z_0)\bs O^-(z_0)$ and all $g\in H(A(z_0)\bs\{z_0\})$. As above, the special structure of the function $\Phi$ will play an important role here. For $z,w\in A(z_0)$, we write $z\sim w$ if there exists some $N\in\mathbb{N}_0$ with $f^N(z)=f^N(w)$. Then $\sim$ defines an equivalence relation on $A(z_0)$, and according to the considerations at the beginning of this section, we have $z\sim w$ if and only if \,$\Phi(z)=\Phi(w)$ (observe that $\Phi|_U=\varphi$ is injective).
\begin{lem}\label{attracting_case_second_inlcusion}
For $g\in H(A(z_0)\bs\{z_0\})$ and $D\subset A(z_0)\bs O^-(z_0)$ open, we have $\omega(D,g,f)\subset H_\Phi(D)$.
\end{lem}
\textbf{Proof:} 
We write $\sim$ for the restriction of the above equivalence relation to $D$ and $[z]_{\sim}$ for the equivalence class of $z\in D$. Moreover, let $D/_{\sim}$ be the quotient induced by $\sim$ and let
$p:D\rightarrow D/_{\sim},\,\,\,p(z):=[z]_{\sim},$
be the associated quotient map. For $h\in\omega(D,g,f)$, there exists a strictly increasing sequence $(n_k)$ in $\mathbb{N}$ such that $(g\circ f^{n_k})$ converges to $h$ locally uniformly on $D$. Hence, $h$ is holomorphic on $D$, and it follows for all points $z,w\in D$ with $z\sim w$ that
$$h(z)=\lim_{k\rightarrow\infty}g(f^{n_k}(z))=\lim_{k\rightarrow\infty}g(f^{n_k}(w))=h(w).$$
Thus, the map $\wt{h}:D/_{\sim}\rightarrow\mathbb{C},\,\wt{h}\big([z]_{\sim}\big):=h(z)$, is well-defined and fulfils $h=\wt{h}\circ p$. Due to the above considerations, the map
$\wt{\Phi}:D/_{\sim}\rightarrow\Phi(D),\,\wt{\Phi}\big([z]_{\sim}\big):=\Phi(z),$
is well-defined, bijective and fulfils $\Phi|_{D}=\wt{\Phi}\circ p$.
Altogether, we have
$$h=\wt{h}\circ p=\wt{h}\circ\wt{\Phi}^{-1}\circ\wt{\Phi}\circ p=\left(\wt{h}\circ\wt{\Phi}^{-1}\right)\circ\Phi|_{D}.$$
As $h$ is holomorphic on $D$ and as $\Phi|_{D}:D\rightarrow\Phi(D)$ is holomorphic and surjective, it follows that $\wt{h}\circ\wt{\Phi}^{-1}$  is holomorphic on $\Phi(D)$ (see \cite{jung}, Lemma A.5 and cf.\,\cite{rudin}, Exercise 14 on p.\,228) so that we obtain $h\in H_\Phi(D)$.\hfill$\square$\medskip

Combining Corollary \ref{attracting_case_first_inlcusion}, Lemma \ref{existence_lemma} and Lemma \ref{attracting_case_second_inlcusion}, we have proved the following concluding statement:
\begin{theo}\label{attracting_case_global}
$G_0$ is an open dense subset of $A(z_0)\bs O^-(z_0)$ such that we have $\dim_H A(z_0)\bs G_0=1$ and
$$\omega(G_0,g,f)=H_\Phi(G_0)$$
for each locally universal function $g$ for $C_f$. 
\end{theo}

In particular, according to the consideration after Definition \ref{locally_universal}, the identity $\omega(G_0,g,f)=H_\Phi(G_0)$ holds for comeager many functions $g\in H(\mathbb{C}\bs\{z_0\})$. Hence the set of all functions in $H(\mathbb{C}\bs\{z_0\})$ which are $H_\Phi(G_0)$-universal for $C_f$ is a comeager set in $H(\mathbb{C}\bs\{z_0\})$.

\subsection{Leau Domains}
\label{subsec:Leau Domains}

Let $z_0\in\mathbb{C}$ be a neutral fixed point of $f$. For $m:=\min\{n\in\mathbb{N}:f^{(n+1)}(z_0)\neq0\}$, there exist $m$ attracting petals $P_1,\ldots,P_m$ for $f$, and for each $k\in\{1,\ldots,m\}$, there exists a conformal map $\varphi_k$ on $P_k$ such that the equation $\varphi_k\circ f^n=\varphi_k+n$ holds on $P_k$ for all integers $n\in\mathbb{N}$ (cf.\,\,Subsection \ref{neutral_case}). For each $k\in\{1,\ldots,m\}$, we consider the open set $G_k:=\bigcup_{n\in\mathbb{N}_0}f^{-n}(P_k)$, which is completely invariant under $f$. As the Leau domain $D_k$ of $f$ which contains $P_k$ (i.e.\,$D_k$ is the component of $F_f$ containing $P_k$) can be written as $D_k=\bigcup_{n\in\mathbb{N}_0}D_{k,n}$, where $D_{k,n}$ is the component of $f^{-n}(P_k)$ which contains $P_k$ (see e.g.\,\cite{steinmetz}, p.\,76), we obtain that $D_k$ is a subset of $G_k$. Analogously to the situation of Subsection \ref{subsec:Schröder Domains}, the image $\Phi_k(G_k)$ is dense in $\mathbb{C}$, and it is well-known that the conjugation map $\varphi_k$ can be extended holomorphically to $G_k$ and that the extended $\varphi_k$ fulfils the same functional equation as $\varphi_k$ (see e.g.\,\cite{milnor}, Corollary 10.9), i.e.\,\,there exists a holomorphic function $\Phi_k:G_k\rightarrow\mathbb{C}$ with $\Phi_k|_{P_k}=\varphi_k$ such that the equation $\Phi_k\circ f^n=\Phi_k+n$ holds on $G_k$ for all $n\in\mathbb{N}$. Hence, for $z,w\in G_k$ and $N\in\mathbb{N}$ with $f^N(z)=f^N(w)$, we have
$$\Phi_k(z)=\Phi_k(f^N(z))-N=\Phi_k(f^N(w))-N=\Phi_k(w),$$
i.e.\,\,the non-constant holomorphic function $\Phi_k$ assumes the same value at all points in $G_k$ which eventually coincide under the iteration of $f$. Considering the subspace
$$H_{\Phi_k}(G_k):=\big\{\psi\circ\Phi_k:\psi\in H(\Phi_k(G_k))\big\},$$
the same approach as in Subsection \ref{subsec:Schröder Domains} yields the following result (cf.\,\cite{jung}, Theorem 4.2.3):
\begin{theo}\label{neutral_case_global_1}
For comeager many functions $g\in H(\mathbb{C}\bs\{z_0\})$, we have $$\omega(G_k,g,f)=H_{\Phi_k}(G_k).$$
\end{theo}

In contrast to the situation of Theorem \ref{attracting_case_global}, we do not have to restrict ourselves to a ``large'' open subset of $G_k$ in order to obtain the identity as stated in Theorem \ref{neutral_case_global_1}. The reason for this is that 
each attracting petal $P_k$ has no holes and is invariant under $f$. Therefore, the inclusion $\omega(G_k,g,f)\supset H_{\Phi_k}(G_k)$ can be proved in exactly the same way as in the proofs of Lemma \ref{attracting_case_generalized} and Corollary \ref{attracting_case_first_inlcusion}.\medskip

Now, let us assume that we have $m\geq2$. Defining $P:=\bigcup_{k=1}^mP_k$ to be the union of all attracting petals, we now consider the set
$$G:=\bigcup_{k=1}^m G_k\,=\,\bigcup_{k=1}^m\bigcup_{n\in\mathbb{N}_0}f^{-n}(P_k)\,=\,\bigcup_{n\in\mathbb{N}_0}f^{-n}(P).$$
As all attracting petals are pairwise disjoint and invariant under $f$, the sets $G_k$ are pairwise disjoint. Moreover, the complete invariance of each set $G_k$ under $f$ implies that $G$ is also completely invariant under $f$. The question arises to determine the sets $\omega(G,g,f)$ for suitable functions $g\in H(G)$. In order to do so, it seems natural to consider the function
$$\Phi:G\rightarrow\mathbb{C},\,\,\,\Phi(z):=\Phi_k(z),\,\,\mbox{if }z\in G_k.$$
As each function $\Phi_k$ quasiconjugates $f|_{G_k}$ to the translation $w\mapsto w+1$ on $\mathbb{C}$ and as each set $G_k$ is invariant under $f$, we obtain that $\Phi$ quasiconjugates $f|_G$ to $w\mapsto w+1$ on $\mathbb{C}$, i.e.\,\,the equation $\Phi\circ f^n=\Phi+n$ holds on $G$ for all integers $n\in\mathbb{N}$. Hence, considering the subspace
$$H_{\Phi}(G):=\big\{\psi\circ\Phi|_D:\psi\in H(\Phi(G))\big\},$$
the following statement can be proved analogously to the first part of the proof of Theorem 4.2.3 in \cite{jung} (using exactly the same techniques as in the proofs of Lemma \ref{attracting_case_generalized} and Corollary \ref{attracting_case_first_inlcusion}):
\begin{lem}\label{neutral_case_sackgasse}
For comeager many functions $g\in H(\mathbb{C}\bs\{z_0\})$, we have $$\omega(G,g,f)\supset H_{\Phi}(G).$$
\end{lem}

However, in this situation, we cannot prove the reverse inclusion as we have done before. The reason for this is the following: Considering on $G$ the equivalence relation
$$z\sim w\,\,\,:\Leftrightarrow\,\,\,\mbox{there exists some }N\in\mathbb{N}_0\mbox{ with }f^N(z)=f^N(w),$$
we obtain that the equality $\Phi(z)=\Phi(w)$ for $z,w\in G$ now does not imply $z\sim w$ anymore (however, due to the above quasiconjugation, $z\sim w$ still yields $\Phi(z)=\Phi(w)$). 
Indeed, let $k,l\in\{1,\ldots,m\}$ with $k\neq l$. As stated above, the intersection $\Phi_k(G_k)\cap\Phi_l(G_l)$ is dense in $\mathbb{C}$ and in particular non-empty. Thus, there exist points $z_k\in G_k$ and $z_l\in G_l$ with $\Phi_k(z_k)=\Phi_l(z_l)$, i.e.\,\,we have $\Phi(z_k)=\Phi(z_l)$. But as $G_k$ and $G_l$ are disjoint and invariant under $f$, we have $f^n(z_k)\neq f^n(z_l)$ for all $n\in\mathbb{N}$ and hence $z_k\not\sim z_l$. 
According to this observation, it might be the case that the set $\omega(G,g,f)$ is a proper superset of the subspace $H_{\Phi}(G)$. In order to avoid the trouble that $\Phi(z)=\Phi(w)$ does not imply $z\sim w$, it is reasonable to consider the function
$$\Phi_{\ast}:G\rightarrow\mathbb{C}\times\{1,\ldots,m\},\,\,\,\Phi_{\ast}(z):=(\Phi_k(z),k),\,\,\mbox{if }z\in G_k.$$
As each function $\Phi_k$ quasiconjugates $f|_{G_k}$ to the translation $w\mapsto w+1$ and as each set $G_k$ is invariant under $f$, we obtain that $\Phi_{\ast}$ quasiconjugates $f|_G$ to the map $(w,k)\mapsto(w,k)+(1,0)$, i.e.\,\,the equation
$\Phi_{\ast}\circ f^n=\Phi_{\ast}+(n,0)$
holds on $G$ for all integers $n\in\mathbb{N}$. Thus, for $z,w\in G$ with $z\sim w$, it follows that $\Phi_{\ast}(z)=\Phi_{\ast}(w)$. Moreover, in contrast to the above considerations about the function $\Phi$, now the reverse implication also holds. Indeed, for $z,w\in G$ with $z\in G_k$, $w\in G_l$ and $\Phi_{\ast}(z)=\Phi_{\ast}(w)$, we have $(\Phi_k(z),k)=(\Phi_l(w),l)$. Hence, we obtain $k=l$ and thus $\Phi_k(z)=\Phi_k(w)$. Analogously to the corresponding considerations in Subsection \ref{subsec:Schröder Domains}, this yields $z\sim w$. Defining
$$H_{\Phi_{\ast}}(G):=\left\{\psi\circ\Phi_{\ast}|_{G}:\psi\in H_{\ast}(\Phi_{\ast}(G))\right\}\vspace{-5pt},$$
where
$$H_{\ast}(\Phi_{\ast}(G)):=\bigcap\limits_{k=1}^m\left\{\psi:\Phi_{\ast}(G)\rightarrow\mathbb{C}: \psi\left(\bull\,|_{\Phi_k(G_k)},k\right)\in H\big(\Phi_k(G_k)\big)\right\},$$
one can show that $H_{\ast}(\Phi_{\ast}(G))$ is a subspace of $H(G)$ which properly contains $H_\Phi(G)$ (cf.\,\cite{jung}, Lemma 4.2.8). Using similar techniques as in Subsection \ref{subsec:Schröder Domains}, the following concluding statement can be proved (cf.\,\cite{jung}, p.\,79ff.):
\begin{theo}\label{neutral_case_global_2}
For comeager many functions $g\in H(\mathbb{C}\bs\{z_0\})$, we have $$\omega(G,g,f)=H_{\Phi_{\ast}}(G).$$
\end{theo}

\subsection{Other Fatou Components}
\label{subsec:rest}

In this final subsection, we want to state several results concerning the question which level of universality the composition operator $C_f$ can have on B\"ottcher domains, Siegel discs, Baker domains and wandering domains of $f$. It is easy to see that the set $\omega\big(G\bs O^-(z_0),g,f\big)$ only contains constant functions if $G$ is a B\"ottcher domain of $f$ which contains a superattracting fixed point $z_0$ of $f$ and $g$ is holomorphic on $\mathbb{C}\bs\{z_0\}$, and that there exist constant functions which are not contained in the closure of $\left\{g\circ f^n|_G:n\in\mathbb{N}\right\}$ in $H(G)$ if $G$ is an invariant Siegel disc of $f$ and $g$ is holomorphic on $G$ (cf.\,\cite{jung}, Remark 4.3.2 and Lemma 4.3.3). Indeed, in these two situations, $C_f$ cannot be $H\big(G\bs O^-(z_0)\big)$-universal or $H(G)$-universal, respectively, because of lack of injectivity of the symbol $f$ on $G\bs O^-(z_0)$ in the first case and because of lack of ``run-away behaviour'' of the sequence of iterates $(f^n)$ on $G$ in the latter case (cf.\,\,Theorem \ref{theorem_grosse_erdmann_mortini}). 
Requiring that $f$ is a transcendental entire function which is injective on a Baker domain or a wandering domain, we again obtain universality of $C_f$. More precisely, the following statements hold (we write $O^+(G):=\bigcup_{k\in\mathbb{N}}f^k(G)$ for the \textit{forward orbit} of $G$ under $f$):
\begin{theo}\quad
\begin{itemize}
\itemsep0pt
\item[i)]Let $G$ be an invariant Baker domain of $f$ and let $f$ be injective on $G$. Then the set of all entire functions which are $H(G)$-universal for $C_f$ is a comeager set in $H(\mathbb{C})$.
\item[ii)]Let $G$ be a wandering domain of $f$, let $f^n$ be injective on $G$ for all $n\in\mathbb{N}$ and let $f^n|_G\rightarrow\infty$. Then the set of all entire functions which are $H(W)$-universal for $C_f$ for all $W\in\mathcal{U}_0(O^+(G))$ is a comeager set in $H(\mathbb{C})$. If, in addition, $G$ is simply connected, the set of all entire functions which are $H(O^+(G))$-universal for $C_f$ is a comeager set in $H(\mathbb{C})$.
\end{itemize}
\end{theo}
\textbf{Proof:}
\begin{itemize}
\itemsep0pt
\item[i)]By definition, $G$ is unbounded and hence simply connected (see \cite{baker}, Theorem 1). Considering the domain $\Omega:=\mathbb{C}$ and the open set $D:=G\subset\Omega$ and observing that the map $f|_D:D\rightarrow D$ is injective and that we have locally uniform convergence $f^n|_D\rightarrow\infty\in\partial_\infty\Omega$ due to the Vitali-Porter theorem, the assertion follows from Corollary \ref{central_universality_theorem_all_open_special_simply_connected}.
\item[ii)]Considering the domain $\Omega:=\mathbb{C}$ and the open set $D:=O^+(G)\subset\Omega$, we obtain that the map $f|_D:D\rightarrow D$ is injective (observe that, as $G$ is a wandering domain of $f$, the sets $f^k(G)$ are pairwise disjoint). 
As we have locally uniform convergence $f^n|_D\rightarrow\infty\in\partial_\infty\Omega$ (cf.\,\cite{jung}, proof of Theorem 4.5.3), the first assertion follows from Corollary \ref{central_universality_theorem_all_open_special}. If $G$ is simply connected, the invariance of the number of holes yields that each set $f^k(G)$ also has no holes (observe that each iterate $f^k$ is injective on $G$). Hence, the second assertion follows directly from the first one.\hfill$\square$
\end{itemize}\smallskip

We remark that there exist several examples of transcendental entire functions which are injective on Baker domains (see e.g.\,\cite{bergweiler_invariant}, p.\,526, \cite{baranski-fagella}, subsection 5.2 and \cite{sienra}, Theorem 1) and diverse examples of transcendental entire functions which have a simply connected wandering domain on which each iterate is injective and on which the iterates converge to the point at infinity (see e.g.\,\cite{eremenko-ljubich}, Example 2 and \cite{zheng}, Theorem 3.4).



\begin{thebibliography}{12}
\parskip0.0mm
\small
\bibitem{baker}
I. N. Baker, \textit{The domains of normality of an entire function}, Ann. Acad. Sci. Fenn. Ser. A I Math. \textbf{1} (1975), 277-283. 
\bibitem{baranski-fagella}
K. Bara\'{n}ski and N. Fagella, \textit{Univalent Baker domains}, Nonlinearity \textbf{14} (2001), 411-429.
\bibitem{bmm}
P. Beise, T. Meyrath and J. M\"uller, \textit{Universality properties of Taylor series inside the domain of holomorphy}, J. Math. Anal. Appl. \textbf{383} (2011), 234-238. 
\bibitem{bergweiler_invariant}
W. Bergweiler, \textit{Invariant domains and singularities}, Math. Proc. Cambridge Philos. Soc. \textbf{117} (1995), 525-532.
\bibitem{bergweiler_uebersicht}
W. Bergweiler, \textit{Iteration of meromorphic functions}, Bull. Amer. Math. Soc. (N.S.) \textbf{29} (1993), 151-188. 
\bibitem{bernal-montes}
L. Bernal-Gonz\'{a}lez and A. Montes-Rodr\'{i}guez, \textit{Universal functions for composition operators}, Complex Variables Theory Appl. \textbf{27} (1995), 47-56.
\bibitem{carleson-gamelin}
L. Carleson and T. W. Gamelin, \textit{Complex Dynamics}, Springer, New York, 1993.
\bibitem{edgar}
G. A. Edgar, \textit{Measure, Topology, and Fractal Geometry}, Springer, New York, 1990.
\bibitem{eremenko-ljubich}
A. E. Eremenko and M. Yu. Lyubich, \textit{Examples of entire functions with pathological dynamics}, J. London Math. Soc. (2) \textbf{36} (1987), 458-468.
\bibitem{grosse-erdmann-diss}
K.-G. Grosse-Erdmann, \textit{Holomorphe Monster und universelle Funktionen}, Mitt. Math. Sem. Giessen \textbf{176} (1987).
\bibitem{grosse-erdmann-survey}
K.-G. Grosse-Erdmann, \textit{Universal families and hypercyclic operators}, Bull. Amer. Math. Soc. (N.S.) \textbf{36} (1999), 345-381.
\bibitem{grosse-erdmann-mortini}
K.-G. Grosse-Erdmann and R. Mortini, \textit{Universal functions for composition operators with non-automorphic symbol}, J. d'Analyse Math. \textbf{107} (2009), 355-376.
\bibitem{jung}
A. Jung, \textit{Universality of Composition Operators with Applications to Complex Dynamics}, Dissertation, University of Trier, URL: http://ubt.opus.hbz-nrw.de/volltexte/2015/963/ (2015).
\bibitem{kahane}
J.-P. Kahane, \textit{Baire's category theorem and trigonometric series}, J. Anal. Math. \textbf{80} (2000), 143-182.
\bibitem{luh}
W. Luh, \textit{Universal functions and conformal mappings}, Serdica \textbf{19} (1993), 161-166.
\bibitem{milnor}
J. Milnor, \textit{Dynamics in One Complex Variable}, Friedr. Vieweg \& Sohn, Braunschweig-Wiesbaden, 1999.
\bibitem{mueller}
J. M\"uller, \textit{Continuous functions with universally divergent Fourier series on small subsets of the circle}, C. R. Math. Acad. Sci. Paris \textbf{348} (2010), 1155-1158.
\bibitem{remmert2}
R. Remmert, \textit{Funktionentheorie 2}, 2., korrigierte Auflage, Springer, Berlin, 1995.  
\bibitem{rudin}
W. Rudin, \textit{Real and Complex Analysis}, Third Edition, McGraw-Hill, New York, 1987.
\bibitem{schiff}
J. L. Schiff, \textit{Normal Families}, Springer, New York, 1993.
\bibitem{schleicher}
D. Schleicher, \textit{Dynamics of entire functions}, In: \textit{Holomorphic dynamical systems}, 295-339, Lecture Notes in Math. \textbf{1998}, Springer, Berlin, 2010. 
\bibitem{sienra}
G. Sienra, \textit{Surgery and hyperbolic univalent Baker domains}, Nonlinearity \textbf{19} (2006), 959-967.
\bibitem{steinmetz}
N. Steinmetz, \textit{Rational Iteration}, de Gruyter, Berlin, 1993.
\bibitem{sullivan}
D. Sullivan, \textit{Quasiconformal homeomorphisms and dynamics I. Solution of the Fatou-Julia problem on wandering domains}, Ann. of Math. \textbf{122} (1985), 401-418.
\bibitem{zheng}
J. H. Zheng, \textit{Iteration of functions which are meromorphic outside a small set}, Tohoku Math. J. (2) \textbf{57} (2005), 23-43. 
\end{thebibliography}
\end{document}